\input amstex
\documentstyle{amsppt}

\exhyphenpenalty=10000

\def
 \today{
  \relax
   \ifcase
    \month
     \or January\or February\or March\or April\or May\or June
      \or July\or August\or September\or October\or November\or December
       \fi
        \space
         \number
          \day,
           \number
            \year}

\def \Aus             {1}
\def \Brown           {2}
\def \Fuchs           {3}
\def \HoltCohomology  {4}
\def \HoltDimension   {5}
\def \HoltOneEnd      {6}
\def \Jech            {7}
\def \KropTwo         {8}
\def \Krop            {9}
\def \KM             {10}
\def \KT             {11}
\def \Osofsky        {12}
\def \Robinson       {13}
\def \Thomas         {14}
\def \tomDieck       {15}

\def \PHKh {\operatorname{\scriptstyle{\text{\bf H}}}}
\def \PHKl {\operatorname{\scriptstyle{\text{\bf L}}}}

\def \arank {\operatorname{\aleph-rank}}

\def \cd {\operatorname{cd}}

\def \E {\operatorname{E}}
\def \uE {\underline{\operatorname{E}}}
\def \End {\operatorname{End}}
\def \Ext {\operatorname{Ext}}
\def \Hom {\operatorname{Hom}}
\def \fHom {\operatorname{\Cal{H}om}}
\def \H {\operatorname{H}}
 \def \C {\operatorname{C}}
  \def \Z {\operatorname{Z}}
   \def \B {\operatorname{B}}
    \def \P {\operatorname{P}}
\def \rank {\operatorname{rank}}
\def \supp {\operatorname{supp}}

\magnification=1200
\vsize 23.7truecm
\vcorrection{-2truecm}
\raggedbottom
\nologo

\topmatter
\leftheadtext{Dicks, Kropholler, Leary  and Thomas}
\rightheadtext{Classifying spaces for proper actions of locally-finite groups}
\title
 Classifying spaces for proper actions of locally-finite groups
  \endtitle
\author
  Warren Dicks, Peter H\. Kropholler, Ian J\. Leary  and Simon Thomas
  \endauthor
\date
 September 5, 2001
  \enddate
\address
  Departament de Matem\`atiques,
   Universitat Aut\`onoma de Bar\-ce\-lo\-na,
    E-08193~Bel\-la\-terra (Bar\-ce\-lo\-na), Spain
     \endaddress
      \email
       dicks\@mat.uab.es
        \endemail
\address
  School of Mathematical Sciences,
   Queen Mary, University of London,
    Mile End Road, London E1 4NS, United Kingdom
     \endaddress
      \email
       P.H.Kropholler\@qmw.ac.uk
        \endemail
\address
  Faculty of Mathematical Studies,
   University of Southampton,
    Southampton SO17 1BJ, United Kingdom
     \endaddress
      \email
       I.J.Leary\@maths.soton.ac.uk
        \endemail
\address
  Department of Mathematics,
   Rutgers University,
    New Brunswick, N.~J.~08903, USA
     \endaddress
      \email
       sthomas\@math.rutgers.edu
        \endemail


\abstract
\nofrills
For each finite ordinal $n$, and each locally-finite group $G$ of cardinality
$\aleph_n$, we construct an $(n+1)$-dimensional, contractible CW-complex on
which
$G$ acts with finite stabilizers.  We use the complex to obtain information
about
cohomology with induced coefficients.  Our techniques also give information
about the
location of some large free abelian groups in the hierarchy~$\PHKh\frak{F}$.
\endabstract

\endtopmatter

\document

Throughout, let $G$ be a group, and let $A$ be a $\Bbb Z$-module with trivial
$G$-action.   We let $AG$ denote the induced $\Bbb ZG$-module
$\Bbb ZG \otimes_{\Bbb Z} A$.  Ring-actions on modules and group-actions on
sets are
tacitly understood to be on the left, if not specified otherwise.

\head 1. Holt's Conjectures \endhead

The purpose of this section is to describe our main algebraic results,
Theorems~3.10
and~5.4, and place them in the context of prior work.

\definition
{1.1 Notation} Let $\rank(G)$ denote the smallest cardinal $\kappa$ such
that there
exists some set of generators of $G$ of cardinality $\kappa$.

If $G$ is {\it not} finitely generated, then $\rank(G) = \vert G\vert$, and
we define
$\arank(G)$ to be the ordinal $\alpha$ such that $\rank(G) = \aleph_\alpha$; if
$G$ {\it is} finitely generated, then $\rank(G) < \vert G\vert$, and we set
$\arank(G)=-1$.

Recall that, for each ordinal $\alpha$,  $\omega_\alpha$ denotes the least
ordinal
of cardinality $\aleph_\alpha$.

We find it convenient to set $\aleph_{-1} = 1$.

Throughout this section,  let $n \in \Bbb N \,(= \omega_0)$. \qed \enddefinition

D\. F\. Holt proposed the following description of the cohomology with induced
coefficients, for locally-finite groups.

\proclaim
{1.2 Conjecture {\rm (Holt~\cite{\HoltOneEnd})} } If $G$ is locally finite,
then
$$\align
\vert \H^n(G,AG)\vert &=  \vert A \vert^{\aleph_{n-1}}
     \;\text{ if }\;  n = \arank(G)+1,\\
\H^n(G,AG) &= 0 \qquad  \quad\text{ if }\; n \ne \arank(G)+1.\endalign$$
\endproclaim

\demo
{Commentary}
In Examples~3.3 we recall that, for any group $G$,
$$\H^n(G,AG) = \cases
A  &\text{ if $G$ is finite and $n=0$},\\
0 &\text{ if $G$ is finite or $n=0$, but not both.}
\endcases\tag"(1.3)"$$
Thus the conjecture really concerns the cases where $n \ge 1$ and $G$ is
infinite,
and the notation has been artificially contrived to embrace the trivial marginal
cases.

For any infinite group $G$, the set of cocycles for $G$ with coefficients in
$AG$ is of cardinality $\vert A\vert^{\vert G\vert}$, so
$1 \le \vert \H^n(G,AG)\vert \le \vert A\vert^{\vert G\vert}$.  If $G$ is
infinite
and locally finite, then the conjecture implies that only the extreme
values can be
achieved.
\qed \enddemo

We now briefly state the cases which are known, including those obtained in this
paper.

\definition
{1.4 Notation}  We say that $G$ has the {\it finite extension property for
proper
subgroups} if each proper subgroup of $G$ is a proper subgroup of finite
index in
some subgroup of $G$.  For example, abelian torsion groups have this
property.

Let us say that $A$ is {\it $o(G)$-inverting} if, for every finite
subgroup $H$ of $G$, multiplication by $\vert H \vert$ gives an
automorphism of $A$;
equivalently, for each $g \in G$ whose order $o(g)$ is finite, multiplication by
$o(g)$ gives an automorphism of $A$.

If $R$ is a ring (associative, with 1), then $R$ is $o(G)$-inverting, as
$\Bbb Z$-module, if and only if the order of each finite subgroup of $G$ is
a unit in
$R$.  If $R$ is not $o(G)$-in\-vert\-ing, then it is easy to show that
$\cd_RG$, the
cohomological dimension of $G$ with respect to $R$, is $\infty$, a value
which we
shall think of as $\omega_0 = \aleph_0$.
\qed \enddefinition

\proclaim
{1.5 Known cases of Conjecture~1.2} Let $G$ be a locally-finite group.
\roster
\item $\H^{n}(G,AG) = 0$ if $n > \arank(G)+1$.
\item If $G$ has the finite extension property for proper subgroups, then
\newline
$\H^n(G,AG) = 0$ if  $n \ne \arank(G) + 1$.
\item For $n\in\{0,1\}$,  $\H^n(G,AG)=0$ if $n \ne \arank(G)+1$.
\item For $n\in\{0,1,2\}$,
$\vert \H^n(G,AG)\vert = \vert A \vert^{\aleph_{n-1}}$ if $n = \arank(G)+1$.
\item   It is consistent with {\rm ZFC} that
\newline
$\vert \H^n(G,AG)\vert \ge  2^{\aleph_{n-1}}$
if $n = \arank(G)+1$ and $A$ is nonzero.
\newline
Hence, it is consistent with {\rm ZFC} that
\newline
$\vert \H^n(G,AG)\vert = \vert A \vert^{\aleph_{n-1}}$ if $n = \arank(G)+1$ and
$\vert A \vert \le \aleph_{n-1}$.
\endroster
\endproclaim

\demo{Commentary} (1), the ``easy" part of Conjecture~1.2, is proved in
Theorem~3.10.  It was proved
in~\cite{\KT} for the case where $A$ is $o(G)$-inverting, and, before that,
in~\cite{\HoltCohomology},~\cite{\HoltDimension} for the case where $A$ is
$o(G)$-inverting and torsion.

(2) was proved by Holt~\cite{\HoltDimension}.  We give another proof of the
abelian
case in Corollary~6.10.

(3). By (1.3), this holds for $n=0$.  It was proved by
Holt~\cite{\HoltOneEnd} for $n=1$; see~Theorem~6.4.

(4).  By (1.3), this holds for $n=0$.  It is well known for $n=1$; see
Theorem~4.5.
In Theorem~5.4, we prove it for $n =2$;  Holt~\cite{\HoltOneEnd} had
previously shown
this was consistent with ZFC, see~\cite{\Thomas,~Section~1}.

(5). Suppose that $n = \arank(G)+1$ and that $A$ is nonzero.

We shall now see that it is consistent with ZFC that
$\vert \H^n(G,AG)\vert \ge  2^{\aleph_{n-1}}$.

For each prime $p$, we write $\Bbb Z(p^\infty):=
\lim\limits_{m\to \infty}\Bbb Z/p^m\Bbb Z$, where, for $m \in \Bbb N$,
the map $\Bbb Z/p^m\Bbb Z \to \Bbb Z/p^{m+1}\Bbb Z$ is given by
multiplication by $p$.

We claim that there exists a $\Bbb Z$-module $k$, and a $\Bbb Z$-submodule $A'$
of $A$, such that the quotient $A/A'$ is isomorphic to $k$, and either
$k = \Bbb Q$, or there exists a prime $p$ such that $k=\Bbb Z/p\Bbb Z$ or
$k= \Bbb Z(p^\infty)$.

Consider first the case where $A$ is {\it not} divisible, so there exists a
prime $p$ such that $A/pA$ is nonzero.  But $A/pA$ is a direct sum of
$\Bbb Z$-submodules each of which is isomorphic to $\Bbb Z/p\Bbb Z$.  Hence
$A/pA$ projects onto any such summand.

If $A$ {\it is} divisible, then $A$ is a  direct sum of $\Bbb Z$-submodules
each of
which is isomorphic to $\Bbb Q$ or to $\Bbb Z(p^\infty)$ for some prime $p$;
see, for example,~\cite{\Fuchs,~Theorem~IV.23.1}.  Hence $A$ projects onto
any such
summand.

In all cases, we can find $A'$, $k$ as claimed.

Now there is a long exact sequence in cohomology which contains the subsequence
$$\H^n(G,AG) \to \H^n(G,kG) \to \H^{n+1}(G,A'G).$$
By (1), $\H^{n+1}(G,A'G) = 0$, so
$\vert \H^n(G,AG) \vert \ge \vert \H^n(G,kG) \vert$.

Thus, for the first part, it remains to show that it is consistent with ZFC
that
$\vert \H^n(G,kG) \vert \ge 2^{\aleph_{n-1}}$.

It is proved in~\cite{\KT} that, if $k$ is an $o(G)$-inverting prime field,
then it
is consistent with ZFC that $\vert\H^n(G,kG)\vert \ge  2^{\aleph_{n-1}}$.
However on carefully reading that proof, one sees that all applications of the
$o(G)$-inverting hypothesis can be replaced with applications of (1), so,
in fact, it
is proved that it is consistent with ZFC that if $k$ is a prime field  then
$\vert\H^n(G,kG)\vert \ge  2^{\aleph_{n-1}}$.

It remains to consider the case where $k = \Bbb Z(p^\infty)$ for some prime
$p$.  Again, it is not difficult to show that the argument in~\cite{\KT} can
be further modified to cover this case by interpreting $\dim k^m:=m$ for
all $m \in \Bbb N$.  If $M$ is a finitely generated (free)
$\Bbb Z$-submodule of $\Bbb ZG$, and $kM$ denotes the image of the natural
map $k\otimes_{\Bbb Z}M \to k\otimes_{\Bbb Z}\Bbb ZG = kG$, then one can show
$\dim kM = \rank M$, since $k$ is divisible.  Moreover, if $M'$ is a
$\Bbb Z$-submodule of $M$, then  $\dim kM' \le \dim kM$, and, if equality
holds, then
$kM' = kM$, since $k$ is divisible.  Using these observations, one can
verify that
the argument in~\cite{\KT} applies with $k = \Bbb Z(p^\infty)$.

It follows that, in all cases, it is consistent with ZFC that
$\vert\H^n(G,kG)\vert \ge  2^{\aleph_{n-1}}$.

Now suppose that $\vert A \vert \le \aleph_{n-1}$. Hence
$\vert A \vert^{\aleph_{n-1}} =  2^{\aleph_{n-1}}$; see~\cite{\Jech, p.49}
for the
case where $n \ge 1$.  Thus, it is consistent with ZFC that
$\vert \H^n(G,AG)\vert \ge \vert A \vert^{\aleph_{n-1}}$.

We previously observed that
$\vert \H^n(G,AG) \vert \le \vert A \vert^{\aleph_{n-1}}$,
so it is consistent with ZFC that
$\vert \H^n(G,AG)\vert  = \vert A \vert^{\aleph_{n-1}}$.

This proves~(5).  It had
previously been proved by Holt~\cite{\HoltDimension} in the case where $G$
has the
finite extension property for proper subgroups; see~\cite{\Thomas,~Section~1}.
\qed\enddemo

We wish to refine part of Conjecture~1.2.

\proclaim
{1.6 Conjecture} If $G$ is locally finite, and  $n = \arank(G)+1$, then
$\H^n(G,AG)$ contains a $\Bbb Z$-sub\-mod\-ule isomorphic to
$A^{\aleph_{n-1}}$,
and hence
$\vert \H^n(G,AG) \vert = \vert A \vert^{\aleph_{n-1}}$.
\endproclaim

\demo
{Commentary} By (1.3), this holds for $n=0$.  It is probably well
known for $n=1$; see Theorem~4.5.  In Theorem~5.4, we prove it for $n = 2$.
\qed
\enddemo

Conjecture~1.2 was preceded by, and motivated by, an earlier proposal,
concerning the
cohomological dimension of locally-finite groups.

\proclaim
{1.7 Conjecture} If $G$ is a locally-finite group, and $R$ is a nonzero,
$o(G)$-inverting ring, then $\cd_RG = \min\{\arank(G)+1, \infty\}.$
\endproclaim

\demo
{Commentary} Holt~\cite{\HoltCohomology} proposed this conjecture with the
additional
hypothesis that $R$ is a field of prime order, and, in~\cite{\KT}, the
additional hypothesis was weakened to $R$ being commutative.

Notice that $\min\{\arank(G)+1, \infty\}$ can be expressed as
$\inf\{n\in\Bbb N\mid \aleph_{n} > \vert G\vert \},$
where the infimum of the empty set is taken to be $\infty$.

The inequality  $\cd_RG \le \arank(G)+ 1$ follows from a classic result of
Auslander~\cite{\Aus,~Proposition~3}; see~\cite{\Osofsky,~Lemma~3.7} or
Theorem~3.10 below.

Cohomological dimension cannot increase on passing to a subgroup, so we may
assume
that $\arank(G) < \omega_0$, and let $n = \arank(G)+ 1$.  The conjecture now
amounts to the claim that $\H^n(G,M) \ne 0$ for some $RG$-mod\-ule $M$.
Notice that,
on $RG$-modules, $\H^{n+1}(G,-)$ vanishes and (hence) $\H^{n}(G,-)$ is
right exact;
also, $M$ is a quotient of some free $RG$-module.  The conjecture is therefore
equivalent to the claim that $\H^n(G,AG) \ne 0$ for some free $R$-module
$A$.  This
claim is implied by the claim that $\H^n(G,AG) \ne 0$ for some vector space
$A$ over
the prime subfield of some simple quotient ring of $R$.  This means that we
may assume
that $R$ is a prime field.

The foregoing claims are known to be consistent with ZFC~\cite{\KT}, and have
been proved in various cases.  The case where $G$ is abelian was proved by
N\. Chen;
see~\cite{\Osofsky,~Corollary~7.6} or Corollary~6.10 below.  Chen's
result has been extended in two directions.
Osofsky~\cite{\Osofsky,~Corollary~7.5}
settled the case where $G$ is generated by finite groups whose pairwise
products are
subgroups; in particular, $G$ has the finite extension property for proper
subgroups.  Holt~\cite{\HoltDimension} settled the case where $G$ has the finite
extension property for proper subgroups and $R$ is a finite prime field.
\qed\enddemo

\head 2. $G$-complexes  \endhead

In this section we construct finite-dimensional contractible spaces with
locally\discretionary{-}{}{-}fi\-nite groups acting on them.

\definition
{2.1 Definitions} A map $f\colon X_1 \to X_2$ between CW-complexes is
{\it cellular} if it carries the $d$-skeleton of $X_1$ to the $d$-skeleton
of $X_2$
for all $d \in \Bbb N$.

A {\it $G$-CW-complex}, or {\it $G$-complex} for short, is a CW-complex $X$
with a
$G$-action such that each element of $G$ acts continuously on $X$,
permuting the open
cells, and fixing only those cells which it fixes pointwise.  It follows
that $G$
acts cellularly.

If $X$ is a $G$-complex then, for each $H \le G$, the set $X^H$, consisting
of points
fixed by  all of $H$, is a CW-subcomplex of $X$, and the set $X/H$
consisting of the
$H$-orbits is a quotient CW-complex.

Let $\frak{sub}(G)$ denote the set of all subgroups of $G$.  A subset
$\frak{X}$ of
$\frak{sub}(G)$ is a {\it subgroup-closed $G$-family} if each subgroup of
each element
of $\frak X$ belongs to $\frak X$ and, moreover, $\frak X$ is closed under
taking
conjugates by elements of $G$.

If $\frak X$ is a subgroup-closed $G$-family, then by
{\it a space of type $\E(G,\frak X)$} we mean a
$G$-complex $X$ with the properties that,
for each $H \in \frak X$,  $X^H$ is contractible, and for each
$H\in \frak{sub}(G) - \frak X$, $X^H$ is empty.  In this event, $X$
is also said to be a {\it classifying space for $G$-actions with stabilizers in
$\frak X$}.

If $\frak X$ is a class of groups, and $\frak X \cap \frak{sub}(G)$ is
a subgroup-closed $G$-family, then by
{\it a space of type $\E(G,\frak X)$} we mean
a space of type $\E(G,\frak X \cap \frak{sub}(G))$.

We let $\frak F$ denote the class of finite groups.  Notice that
$\frak F \cap \frak{sub}(G)$ is a sub\-group\discretionary{-}{}{-}closed
$G$-family.  A space of type $\E(G,\frak F)$ is called {\it an
$\uE G$}.  (It is also called a classifying space for proper $G$-actions,
that is, $G$-actions with finite stabilizers.)
\qed \enddefinition

The following is well known.

\proclaim
{2.2 Proposition} If $\frak X$ is a subgroup-closed $G$-family,
then there exists a space of type $\E(G,\frak X)$, and any $G$-map
between two spaces of type $\E(G,\frak X)$ is a $G$-homotopy
equivalence.
\endproclaim

\demo{Proof}  The first part can be seen by Milnor's construction. Thus,
let $\Delta$
be any $G$-set such that $\frak X$ is precisely the set of subgroups of $G$
which fix
at least one point of $\Delta$.   Let
$X = \Delta \ast \Delta \ast \Delta \ast \cdots$,
the union of iterated joins of $\Delta$.  Then $X$ is a space of type
$\E(G,\frak X)$.

For the second part, see, for example,~\cite{\tomDieck,~Proposition~II.2.7}.
\qed \enddemo

\proclaim
{2.3 Corollary} If $\frak X_1 \subseteq \frak X_2$ are subgroup-closed
$G$-families,
and
$X_1$ {\rm(}resp. $X_2${\rm)} is a space of type
$\E(G,\frak X_1)$ {\rm(}resp. $\E(G,\frak X_2)${\rm)}, then
there exists a cellular $G$-map $X_1 \to X_2$.
\endproclaim

\demo{Proof} The join $X_1 \ast X_2$ is a space of type
$\E(G,\frak X_2)$, and the inclusions
$$\iota_1\colon X_1 \to X_1 \ast X_2, \quad \iota_2\colon X_2 \to X_1 \ast X_2$$
are $G$-maps.  By Proposition~2.2,
$\iota_2$ is a $G$-homotopy equivalence, and the homotopy inverse
$X_1\ast X_2 \to X_2$ composed with $\iota_1$ gives a $G$-map  $X_1 \to
X_2$. This is
then $G$-homotopic to a cellular $G$-map $X_1 \to X_2$; see, for
example,~\cite{\tomDieck,~Theorem~II.2.1}.
\qed
\enddemo

One could give a dual proof, using the projection maps from the Cartesian
product $X_1 \times X_2$, which is a space of type $\E(G,\frak X_1)$.

The following is a topological analogue of a classic result of
Auslander~\cite{\Aus,~Proposition~3}.

\proclaim
{2.4 Theorem} Let $\beta$ be a limit ordinal, let
$(G_\alpha \mid  \alpha \le \beta)$  be a continuous chain of subgroups of $G$,
and let $(\frak{X}_\alpha \mid  \alpha \le \beta)$ be a continuous chain of
subsets of $\frak{sub}(G)$ such that, for each $\alpha \le \beta$,
$\frak{X}_\alpha$ is a subgroup-closed $G_\alpha$-family.

Let $n \in \Bbb N$, and suppose that, for each $\alpha < \beta$, there exists an
$n$-dimensional space $Y_\alpha$ of type  $\E(G_\alpha, \frak{X}_\alpha)$.
Then there
exists an $(n+1)$-dimensional space of type  $\E(G_\beta,\frak{X}_\beta)$.
\endproclaim

\demo{Proof} For each $\alpha < \beta$,
$\frak X_\alpha \subseteq \frak X_{\alpha+1} \cap \frak{sub}(G_\alpha)$ are
subgroup-closed $G_\alpha$-families.  Also, $Y_\alpha$  is a
space of type $\E(G_\alpha, \frak X_\alpha)$, and  $Y_{\alpha+1}$ can be
viewed as a
space of type $\E(G_\alpha, \frak X_{\alpha+1}\cap \frak{sub}(G_\alpha))$. By
Corollary~2.3, there exists a cellular $G_{\alpha}$-map
$Y_{\alpha} \to Y_{\alpha+1}$,  and hence a cellular $G_{\alpha+1}$-map
$f_{\alpha}\colon G_{\alpha+1}\times_{G_\alpha}Y_{\alpha} \to Y_{\alpha+1}$.
Let $M_{\alpha}$ denote the mapping cylinder of $f_{\alpha}$.  Since
$f_{\alpha}$ is cellular, $M_{\alpha}$ has the structure of a CW-complex,
and is a
space of type $\E(G_{\alpha+1},\frak X_{\alpha+1})$. Notice that $\dim
M_{\alpha}
= n+1$, since $\dim Y_\alpha = \dim Y_{\alpha+1} = n$.

We recursively construct a continuous chain
$(X_\alpha \mid \alpha < \beta)$ where
$X_\alpha$ is a space of type $\E(G_\alpha,\frak X_\alpha)$, and, for
$\alpha \ge 1$,
$\dim X_\alpha= n+1$.

We take $X_0  = Y_0$, and at limit ordinals we take directed unions.

Suppose $\alpha < \beta$ and that $X_\alpha$ has been constructed.

By Proposition~2.2, since $X_\alpha$ and $Y_\alpha$ are of type
$\E(G_\alpha,\frak X_\alpha)$ there exists a cellular $G_\alpha$-map
$Y_\alpha \to X_\alpha$, and hence a  cellular $G_{\alpha+1}$-map
$G_{\alpha+1} \times_{G_\alpha} Y_\alpha \to G_{\alpha+1}
\times_{G_\alpha}X_\alpha$.
Take $X_{\alpha+1}$ to be the identification space, or pushout,
$$ \CD G_{\alpha+1}\times_{G_\alpha} Y_{\alpha}  @>>>  M_{\alpha} \\
@VV  V   @VV V\\
G_{\alpha+1}\times_{G_\alpha} X_{\alpha} @>>>  X_{\alpha+1}.
\endCD$$
Notice that $\dim X_{\alpha+1}= n+1$, since $\dim Y_{\alpha} =n$, $\dim
M_{\alpha} =
n+1$, and \linebreak $\dim X_\alpha \le n+1$.   It is not difficult to
check that
$X_{\alpha+1}$ is  of type  $\E(G_{\alpha+1},\frak X_{\alpha+1})$.

This completes the proof.
\qed
\enddemo

\definition
{2.5 Remark}  For $n=0$ and $\beta = \omega_0$, the construction in the above
proof gives the Bass-Serre tree of the graph of groups corresponding to the
countable
ascending chain $(G_\alpha \mid  \alpha < \omega_0)$.  Here,
$\frak X_{\omega_0} = \bigcup\limits_{\alpha < \omega_0}
\frak{sub}(G_\alpha)$.
\qed\enddefinition

\proclaim
{2.6 Theorem} If $n \in \Bbb N$, and $G$ is a locally-finite group with
$\arank(G) < n$, then there exists an $n$-dimensional $\uE G$.
\endproclaim

\demo{Proof} We argue by induction on $n$.

If $n=0$, then $G$ is finite.  Here a single point with trivial $G$-action is a
$0$-dimensional $\uE G$.

Thus we may assume that $n \ge 1$, and that the result holds for smaller $n$.

We can choose a continuous chain $(G_\alpha \mid \alpha \le \omega_n)$ of
subgroups of $G$ such that $G_{\omega_n} = G$, and, for each $\alpha <
\omega_n$,
$\arank(G_\alpha) < n-1$, so, by the induction hypothesis, there exists an
$(n-1)$-dimensional $\uE G_\alpha$.  Thus, by Theorem~2.4, there
exists an $n$-dimensional  $\uE G$.

This completes the proof.
\qed
\enddemo

\definition
{2.7 Remarks}  The foregoing construction applies in greater generality.

Suppose that, for every group $H$, there is specified a subgroup-closed
$H$-family
$\frak Y(H)$ satisfying the following three conditions:

\itemitem{} Any group isomorphism $H_1 \to H_2$ induces a bijection
$\frak Y(H_1) \to \frak Y(H_2)$.

\itemitem{} If $H_1 \le H_2$, then $\frak Y(H_1) \subseteq \frak Y(H_2)$.

\itemitem{} If $H$ is the union of a well-ordered chain of subgroups
$H_\alpha$, then
$\frak Y(H)$ is the union of the  $\frak Y(H_\alpha)$.

Let $\frak G_0$ denote the class consisting of those groups $H$ such that
$\frak Y(H) = \frak{sub}(H)$.  For $n < \omega_0$, recursively define
$\frak G_{n+1}$ to be the class consisting of those groups which can be
expressed as
the union of a well-ordered chain of subgroups which lie in $\frak G_{n}$.

The above argument then shows that if $G \in \frak G_{n}$, then there exists an
$n$-di\-men\-sion\-al $\E(G, \frak Y(G))$.

For any $n > 0$ and any $G \in \frak G_{n}$, it can be arranged that all of the
spaces involved in the construction of $\E(G, \frak Y(G))$ have distinguished
contractible subcomplexes which are transversals for the group actions.
In particular, the quotient complex $\E(G, \frak Y(G))/G$ is contractible,
although we will not use this information.
\qed \enddefinition

We record one example.

\proclaim
{2.8 Theorem}
If $n \in \Bbb N$, and $\arank(G) < n$, and $\frak X$ is the set of all
subgroups of
all finitely generated subgroups of $G$, then there exists an $n$-dimensional
$\E(G,\frak X)$.
\endproclaim

\demo{Proof} For each group $H$, let $\frak Y(H)$ be the set consisting of the
subgroups of the finitely generated subgroups of $H$.  It is easy to see that
$\frak Y(-)$ respects isomorphisms, inclusions and well-ordered unions.

It can be shown, by induction on $n$, that $G \in \frak G_n$, in
the notation of the previous remark, so, by that remark, there exists an
$n$-dimensional $\E(G,\frak Y(G))$.
\qed \enddemo

\definition
{2.9 Example}  In the foregoing theorem, if $G$ is abelian (or locally
finite), then
$\frak X$ is the set of finitely generated subgroups of $G$.
\qed \enddefinition

\head 3. Eventual vanishing of cohomology with induced coefficients \endhead

In this section, we recall how $\H^*(G,-)$ can be computed using an
$\uE G$, and apply the method in the case where $G$ is locally finite.

\definition
{3.1 Definitions}  Let $M$ be a $\Bbb ZG$-module.

We say that $M$ is {\it $G$-acyclic} if $\H^n(G,M) = 0$ for all $n \ge 1$.

Any $\Bbb ZG$-summand of a $G$-acyclic $\Bbb ZG$-module is again $G$-acyclic.

If $\quad\cdots \to P_1 \to P_0 \to \Bbb Z \to 0$ is a $\Bbb ZG$-projective
resolution
of
$\Bbb Z$, then it is easy to see that $M$ is $G$-acyclic if and only if the
sequence
$$0 \to \Hom_{\Bbb ZG}(\Bbb Z,M) \to \Hom_{\Bbb ZG}(P_0, M) \to
\Hom_{\Bbb ZG}(P_1, M) \to \cdots $$ is exact.   When this holds, we say that
$\Hom_{\Bbb ZG}(-,M)$ {\it carries augmented $\Bbb ZG$-pro\-jec\-tive
resolutions of
$\Bbb Z$ to exact sequences}.

For any set $\Delta$, we let $A[[\Delta]]$ denote the set of all functions from
$\Delta$ to $A$, and such a function $x \mapsto a_x$ will be written as a formal
sum $\sum\limits_{x\in \Delta}a_x.x$, and its {\it support} is
$$\supp(\sum\limits_{x\in \Delta}a_x.x)\colon= \{x \in \Delta \mid a_x \ne
0\}.$$
The $\Bbb Z$-module structure of $A[[\Delta]]$ is defined in the obvious
way.  We let
$A[\Delta]$, or $A\Delta$ denote the  $\Bbb Z$-submodule of $A[[\Delta]]$
consisting
of the elements with finite support. If $\Delta$ is a $G$-set, we view
$A[[\Delta]]$ as a $\Bbb ZG$-module, with
$$g(\sum\limits_{x\in \Delta} a_x.x) = \sum\limits_{x\in \Delta} a_x.gx
= \sum\limits_{x\in \Delta} a_{g^{-1}x}.x.$$
(Notice that the action of $G$ on the set of all functions from $\Delta$ to
$A$ is
traditionally, and more naturally, on the right, with
$$(x \mapsto a_x )g =  (x \mapsto a_{gx}),$$
and our action of $g$ on the left corresponds to the usual action of
$g^{-1}$ on the
right.  For our purposes, it is convenient to have actions on the left, wherever
possible.)  Then $A\Delta$ is a $\Bbb ZG$-submodule of $A[[\Delta]]$.  In
the case
where $\Delta = G$, this notation is consistent with the notation for the
induced
$\Bbb ZG$-module $AG$.  Here $A[[G]]$ is called a {\it coinduced}
$\Bbb ZG$-module.

There is a natural bijection
$$\Hom_{\Bbb Z}(M, A) \to \Hom_{\Bbb ZG}(M,A[[G]]),
\quad \psi \mapsto (m \mapsto \sum\limits_{x\in G}
\psi(x^{-1}m).x).\tag"(3.2)"$$
Since $\Hom_{\Bbb Z}(-, A)$ carries $\Bbb Z$-split exact
sequences of $\Bbb ZG$-modules to $\Bbb ZG$-split exact
sequences of $\Bbb ZG$-modules, we see that co-induced
$\Bbb ZG$-modules are $G$-acyclic.
\qed\enddefinition

\definition
{3.3 Examples}  Let $G$ be a finite group.

Here, $AG = A[[G]]$, so induced $\Bbb ZG$-modules are co-induced, and hence
$G$-acyclic.

Suppose that $M$ is an $o(G)$-inverting $\Bbb ZG$-module.  Then the
multiplication map
$$M[[G]] \to M, \quad\sum_{g\in G} m_g.g\mapsto \sum_{g\in G}gm_g,$$
is
$\Bbb ZG$-split with right inverse
$m \mapsto \frac{1}{\vert G\vert}\sum\limits_{g\in G} g^{-1}m.g$.  Here,
$M$ is a
$\Bbb ZG$-summand of an induced $\Bbb ZG$-module, so $M$ is $G$-acyclic.
\qed\enddefinition

In the following, $G$ acts on tensor products over $\Bbb Z$ via the
diagonal action.

\proclaim
{3.4 Lemma} Let $M$ be a $\Bbb ZG$-module.
\roster
\item The functor $\Hom_{\Bbb ZG}(- \otimes_\Bbb Z\Bbb ZG, M)$ carries
$\Bbb Z$-split exact sequences of $\Bbb ZG$-mo\-dules to exact sequences.
\item Let $H$ be a subgroup of $G$.  If $M$ is $H$-acyclic, then the
functor
$$\Hom_{\Bbb ZG}(\Bbb Z[G/H]\otimes_\Bbb Z -, M)$$ carries augmented
$\Bbb ZG$-projective resolutions of $\Bbb Z$ to exact sequences.
\endroster
\endproclaim

\demo{Proof}  (2). Let $L$ be a $\Bbb ZG$-module.  There is a natural
identification of $\Bbb ZG$-modules,
$$\Bbb Z[G/H]\otimes_{\Bbb Z}L = \Bbb ZG\otimes_{\Bbb ZH}L,$$
with $gH\otimes\ell$ corresponding to $g\otimes g^{-1}\ell$.

It follows that we can identify
$$\Hom_{\Bbb ZG}(\Bbb Z[G/H]\otimes_\Bbb Z -, M) =
\Hom_{\Bbb ZG}(\Bbb ZG\otimes_{\Bbb ZH} (-, M)
= \Hom_{\Bbb ZH}(-, M)$$ as functors on $\Bbb ZG$-modules.
Since $M$ is $H$-acyc\-lic, this functor carries augmented $\Bbb
ZG$-projective resolutions of $\Bbb Z$ to exact sequences.

(1) is proved similarly. \qed
\enddemo

\definition
{3.5 Notation} Let $X$ be a $G$-complex.

We shall treat $X$ as the $G$-set whose elements are the open cells of $X$.  The
cellular chain complex of $X$ is then the permutation module $\Bbb ZX$,
with the structure of a differential graded $\Bbb ZG$-module, with differential
$\partial$ of degree $-1$.  Here the grading is that determined by the
dimensions of
the cells, so the $n$th component $\C_n(\Bbb ZX)$ has as $\Bbb Z$-basis the
cells of
dimension $n$.

We let $\eta\colon X\times X\to\Bbb Z$ denote the function such that
$\partial x= \sum\limits_{y \in X} \eta(x,y).y$  for each $x \in X$.  Thus, if $x$ is
an $n$-cell, then $\eta(x,y) = 0$ unless $y$ is one of the finitely many
$(n-1)$-cells incident to $x$.
\qed \enddefinition

The following is a degenerate case of the equivariant cohomology spectral
sequence;
see, for example,~\cite{\Brown, VII.7.10(7.10)}.

\proclaim
{3.6 Theorem} Let $X$ be an acyclic $G$-complex. If $M$ is a $\Bbb
ZG$-module which is
$G_x$-acyclic for each  $x \in X$, then
$\H^*(\fHom_{\Bbb ZG}(\Bbb Z X,M))  \simeq \H^*(G,M)$,  as graded abelian
groups.
\endproclaim

Recall that $\fHom_{\Bbb ZG}(\Bbb Z X,M)$ denotes the differential graded
abelian group with $n$th component
$\C^n(\fHom_{\Bbb ZG}(\Bbb Z X,M)) = \Hom_{\Bbb ZG}(\C_n(\Bbb Z X),M)$.

\demo{Proof}  The homology of $(\Bbb ZX, \partial)$ is $\Bbb Z$, concentrated in
degree zero.

We choose a free $\Bbb ZG$-resolution of $\Bbb Z$, and write it as
$(\Bbb ZY, \partial)$ for some $G$-free $G$-set $Y$; for
example, we could take $Y$ to be an $\E G$, and $\Bbb ZY$ its cellular
chain complex.  Then
$\fHom_{\Bbb ZG}(\Bbb ZY, M)$ is an additive abelian differential graded
group, and
its cohomology is  $\H^*(G,M)$.

We consider the double complex $\Bbb ZX \otimes_{\Bbb Z} \Bbb ZY$ with diagonal
$G$-action, and the double complex
$\fHom_{\Bbb Z G}(\Bbb Z X \otimes_{\Bbb Z} \Bbb ZY, M)$.
We get a fourth-quadrant commuting diagram which can be schematically
represented as
$$ \CD
@. @. 0 \\
@.  @. @VVV\\
@. @. \fHom_{\Bbb Z G}(\Bbb ZX, M)\\
@.  @.  @VVV \\
0 @>>> \fHom_{\Bbb Z G}(\Bbb ZY, M) @>>>
\fHom_{\Bbb Z G}(\Bbb ZX \otimes_{\Bbb Z} \Bbb ZY, M).
\endCD \tag"(3.7)"$$

To show that the cohomology group of the outer row,
$\fHom_{\Bbb Z G}(\Bbb ZX,M)$, is isomorphic to the cohomology group of the
outer column, $\fHom_{\Bbb Z G}(\Bbb ZY, M)$, it suffices
to show that the remaining, or inner, rows and columns of (3.7) are exact.
Each inner
column is exact because
$\fHom_{\Bbb Z G}(\Bbb ZX \otimes_{\Bbb Z} -, M)$ is exact on augmented
projective
$\Bbb ZG$-resolutions of $\Bbb Z$, by Lemma~3.4(2).  Similarly, each inner row
is exact because $\fHom_{\Bbb Z G}(- \otimes_{\Bbb Z} \Bbb ZY, M)$ is exact on
$\Bbb Z$-split exact sequences of $\Bbb ZG$-modules, by Lemma~3.4(1). \qed
\enddemo

\proclaim
{3.8 Corollary} Let $M$ be a $\Bbb ZG$-module, and let $X$  be a
finite-di\-men\-sional acyclic $G$-complex.  If $M$ is $G_x$-acyclic for each
$x \in X$, then  $\H^{n}(G,M) = 0$ for all \linebreak $n > \dim X$.
\qed \endproclaim

We record the case of finite stabilizers.

\proclaim
{3.9 Corollary}  Let $M$ be a $\Bbb ZG$-module, and suppose that $M$ is
$H$-acyclic for each finite subgroup $H$ of $G$; for example, this holds if
$M = AG$,
or if $M$ is
$o(G)$-inverting.  Let $X$ be an acyclic $G$-complex with finite
stabilizers;  for
example, this holds if $X$ is an $\uE G$.  Then $\H^{n}(G,M) = 0$ for all
$n > \dim X$.   \qed \endproclaim

Here we can apply Theorem~2.6.

\proclaim
{3.10 Theorem}  Let $G$ be a locally-finite group, and  $M$ a $\Bbb
ZG$-module which
is $H$-acyclic for each finite subgroup $H$ of $G$; for example, this holds
if $M =
AG$, or if $M$ is  $o(G)$-inverting.  Then $\H^{n}(G,M) = 0$ for all $n >
\arank(G) +
1$.
\qed
\endproclaim

\definition
{3.11 Remark} Theorem~3.10 can also be proved using the argument of the
first paragraph of~\cite{\KT, Section 1}. \qed
\enddefinition

\head 4. Locally-finite groups of cardinality $\aleph_0$  \endhead

In this section, we recall how $\H^*(G,AG)$ can be computed using an
$\uE G$, and apply the method in the one-dimensional case.

\definition
{4.1 Definitions} Let $X$ be an $\uE G$, or, more generally, any acyclic
$G$-complex in which all cell stabilizers are finite, and let Notation~3.5
apply.

We have natural identifications
$$A[[X]] = A^X = \Hom_{\Bbb Z}(\Bbb ZX, A).\tag"(4.2)"$$

For simplicity, let us suppose that $X$ is finite dimensional.

Then $A[[X]]$ has the structure of a differential graded
$\Bbb ZG$-module, in which the differential $\partial^*$ has degree +1, and
is given
by
$$\partial^*(\sum\limits_{x\in X} a_x.x)
= \sum\limits_{x\in X} (\sum\limits_{y\in X} \eta(x,y) a_y).x.$$

The cohomology of $(A[[X]], \partial^*)$ is $A$ concentrated in degree zero.

Let $$A_G[[X]]:= \{\sum\limits_{x\in X} a_x.x \in A[[X]]
 \mid \{g\in G \mid a_{gx} \ne 0\} \text{ is finite, for all }x \in X\}.$$

Since $G$-stabilizers are finite, we see that $A_G[[X]]$ consists of all
functions from $X$ to $A$ with finite support in each $G$-orbit.

It is straightforward to check that $A_G[[X]]$ is a differential graded
$\Bbb ZG$-submodule of $A[[X]]$.

We write $\C^n(A_G[[X]])$, $\B^n(A_G[[X]])$, and $\Z^n(A_G[[X]])$
for the $n$-cochains, $n$-co\-boundaries, and $n$-cocycles, respectively.
\qed\enddefinition

Sometimes the notation $\fHom_c(\Bbb Z[X], A)$ is used to denote
$A_G[[X]]$; see, for
example,~\cite{\Brown,~Lemma~VIII.7.4}.

The following is a variation on the usual ``compact supports" cohomology;
see, for
example,~\cite{\Brown,~Proposition~VIII.7.5}.  It is particularly useful in
the study
of ends of groups.

\proclaim
{4.3 Theorem} If $X$ is a finite-dimensional acyclic $G$-complex with finite
stabilizers, then there is a natural isomorphism $\H^*(A_G[[X]]) \simeq
\H^*(G,AG)$
of graded abelian groups.
\endproclaim

\demo{Proof} There is a natural identification of
$\Hom_{\Bbb ZG}(\Bbb ZX, A[[G]])$ with
 $\Hom_{\Bbb Z}(\Bbb ZX, A)$; see~(3.2).
There is also a natural identification of
$\Hom_{\Bbb Z}(\Bbb ZX, A)$ with $A[[X]]$; see~(4.2).   It is
easy to show that under these identifications, $\fHom_{\Bbb ZG}(\Bbb ZX, AG)$
corresponds to $A_G[[X]]$.  Hence
$\H^*(A_G[[X]]) \simeq \H^*(\fHom_{\Bbb ZG}(\Bbb ZX, AG))$.  Finally,
$\H^*(\fHom_{\Bbb ZG}(\Bbb ZX, AG))\simeq \H^*(G,AG)$, by Theorem~3.6.
\qed\enddemo

There is a natural right $G$-action on $\H^*(G,AG)$, arising from the
$\Bbb ZG$-bimodule structure on $AG$.  This agrees with the natural right
$G$-action
on $A_G[[X]]$ which we have transformed into a left $G$-action.

Let us illustrate how Theorem~4.3 can be used to study $\H^*(G,AG)$ when $G$ is
locally finite of cardinality $\aleph_0$.  To do this we now construct a
standardized $\uE G$, as in Remark~2.5.

\definition
{4.4 Definition}  Let $G$ be a locally-finite group of cardinality $\aleph_0$.

Index the elements of $G$ with $\omega_0$, so
$G = \{h_\alpha \mid \alpha < \omega_0\}.$
For each $\beta \le \omega_0$, let
$H_\beta: = \langle h_\alpha \mid \alpha < \beta\rangle$.
Let $(G_\alpha \mid \alpha \le \omega_0)$ be the subsequence of
$(H_\alpha \mid \alpha \le \omega_0)$ obtained by omitting each term which is
equal to an earlier term.

Notice that $G_{\omega_0}= G$ and, for each $\alpha < \omega_0$,
$\vert G_\alpha \vert < \aleph_0$.  Moreover,  $G_0 = 1$, and, for each
$\alpha < \omega_0$,
$G_{\alpha+1} = \langle G_\alpha, g_\alpha \rangle$, where $g_\alpha =
h_{\alpha'}$
and $\alpha'$ is the least ordinal such that $h_{\alpha'} \notin G_\alpha$.

We define the {\it line}, denoted  $\Bbb R$, to be the tree with
vertices  $(v_n\mid  n \in \Bbb Z)$, (oriented) edges
$(e_n \mid n \in \Bbb Z)$, and incidence relations
$\iota e_n = v_n$, $\tau e_n = v_{n+1}$, for all
$n \in \Bbb Z$.

We define the {\it half-line}, denoted  $\Bbb R^+$, to be the subtree of
$\Bbb R$ with vertices \linebreak $(v_\alpha\mid  \alpha < \omega_0)$, and
edges
$(e_\alpha \mid \alpha < \omega_0)$.

There is a $G$-tree $X$ with $G$-transversal $\Bbb R^+$ such that the
$G$-stabilizer of $v_\alpha$, and of $e_\alpha$, is $G_\alpha$, for all
$\alpha < \omega_0$.  This completely specifies $X$.  It is not difficult
to see that $X$ is a locally-finite tree with one end.

We denote the set of vertices (resp. edges) of $X$ by $VX$ (resp. $EX$).
\qed\enddefinition

The following is fairly standard, but we do not know an explicit reference.

\proclaim
{4.5 Theorem}  If $G$ is a locally-finite group of cardinality $\aleph_0$,  then
$\H^*(G,AG)$ is concentrated in degree $1$,
$\H^1(G,AG)$ contains a $\Bbb Z$-submodule isomorphic to $A^{\aleph_0}$, and
$\vert \H^1(G,AG) \vert
= \vert A\vert^{\aleph_0} = \vert A\vert^{\vert G \vert}$.
\endproclaim

\demo {Proof} By (1.3) and Theorem~3.10, $\H^*(G,AG)$ is concentrated in
degree 1, so
it remains to study $\H^1(G,AG)$.

Let $X$ be as in Definition~4.4.

By Theorem~4.3, $\H^1(G,AG) = \H^1(A_G[[X]])$.

Let $Q$ denote $\{e_\alpha \mid \alpha< \omega_0\}$, so we can view
$$A[[Q]]\subseteq A_G[[EX]] = \C^1(A_G[[X]]) = \Z^1(A_G[[X]]).$$

Let $\Cal P$ be a partition of $\omega_0$ into $\aleph_0$ subsets, each of
cardinality $\aleph_0$.  We denote the map $\omega_0 \to \Cal P$ by $\alpha
\mapsto
[\alpha]$.

Consider any $\phi \in A^{\Cal P}$, $[\alpha] \mapsto \phi[\alpha]$.

Define $$\phi^\dag\colon =
\sum\limits_{\alpha < \omega_0} \phi[\alpha].e_\alpha
\in A[[Q]] \subseteq \Z^1(A_G[[X]]).$$

We claim that if $\phi^\dag \in \B^1(A_G[[X]])$ then $\phi = 0$.

Suppose not, so there exists $\psi \in \C^0(A_G[[X]]) = A_G[[VX]]$ such that
$\partial^*\psi = \phi^\dag$, and there exists $p \in \Cal P$ such that
$\phi(p) \ne 0$.

There exists $\mu < \omega_0$ such that
$\supp(\psi) \cap Gv_0 \subseteq G_\mu v_0$.

There exists $\alpha \in p$ such that $\alpha > \mu$. Notice
$\phi[\alpha] = \phi(p) \ne 0$.

In $\Bbb ZX$, $\partial(\sum\limits_{i=0}^\alpha e_i) = v_{\alpha+1} -
v_0$.  Hence
$\partial((g_{\alpha}-g_{\alpha-1})\sum\limits_{i=0}^\alpha e_i) =
-(g_{\alpha}-g_{\alpha-1})v_0$, since $g_\alpha, g_{\alpha-1}
\in G_{\alpha+1} = G_{v_{\alpha+1}}$.

We can view $\psi$ as an additive map $\Bbb ZX \to A$, and apply it to the
foregoing
to get
$\phi^\dag((g_{\alpha}-g_{\alpha-1})\sum\limits_{i=0}^\alpha e_i)
  = \psi(-(g_{\alpha}-g_{\alpha-1})v_0) = 0$, since
$g_\alpha, g_{\alpha-1} \not\in G_{\alpha-1} \supseteq G_\mu$.
Thus $$0 = \phi^\dag((g_{\alpha}-g_{\alpha-1})\sum\limits_{i=0}^\alpha e_i)
= \phi^\dag(\sum\limits_{i=0}^\alpha g_{\alpha}e_i)
- \phi^\dag(\sum\limits_{i=0}^{\alpha-1} g_{\alpha-1}e_i)
- \phi^\dag(e_\alpha) = 0 - 0 - \phi(p).$$

Hence $\phi(p) = 0$, which is a contradiction.

This proves that the composition
$$A^{\Cal P} \to A[[Q]] \to \H^1(A_G[[X]]) = \H^1(G,AG)$$ is injective.
Since $\vert \Cal P \vert = \aleph_0$,
$A^{\aleph_0}$ embeds in $\H^1(G,AG)$.
\qed \enddemo

\head 5. Locally-finite groups of cardinality $\aleph_1$   \endhead

In this section we study $\H^*(G,AG)$ when $G$ is locally finite with
$\arank(G) = 1$,
topologizing and refining results of D\. J\. Holt.

We begin by constructing a standardized $\uE G$.

\definition
{5.1 Definitions} Let $G$ be a locally-finite group of cardinality $\aleph_1$.

Let $\omega_1'$ denote the set of limit ordinals less than $\omega_1$.

As in Definition~4.4, we start by indexing the elements of $G$,
$G = \{h_\alpha \mid \alpha < \omega_1\}$, set
$H_\beta: = \langle h_\alpha \mid  \alpha < \beta\rangle$ for each
$\beta \le \omega_1$, and let $(G_\alpha \mid \alpha \le \omega_1)$ be the
subsequence of $(H_\alpha \mid \alpha \le \omega_1)$ obtained by omitting
each term
which is  either finite, or equal to an earlier term.

Notice $(G_\alpha \mid \alpha \le \omega_1)$ is a continuous chain of
subgroups of $G$, $G_{\omega_1}= G$, and, for each $\alpha < \omega_1$,
$\vert G_\alpha \vert = \aleph_0$ and
$G_{\alpha+1} = \langle G_\alpha, g_\alpha \rangle$,  where $g_\alpha =
h_{\alpha'}$,
and $\alpha'$ is the least ordinal such that $h_{\alpha'} \notin G_\alpha$.

For each subgroup $H$ of $G$, we set $H^\bullet:= G_0 \cap H$.

We shall now construct a family
$(G_{\alpha,n} \mid  \alpha < \omega_1,  n < \omega_0)$ of finite subgroups
of $G$
such that the following hold:
\roster
\item for each $\alpha < \omega_1$,
$(G_{\alpha,n} \mid  n < \omega_0)$ is an
increasing chain with union $G_\alpha$;
\item for each $\alpha < \omega_1$, $n < \omega_0$,
$G_{\alpha+1,n} = \langle G_{\alpha,n}, g_\alpha \rangle$;
\item for each $\alpha \in \omega_1'$, $n < \omega_0$, there is a distinguished
element
$$g_{\alpha,n} \in G_{\alpha,n+1}^\bullet - G_{\alpha,n}^\bullet
(=G_0 \cap (G_{\alpha,n+1}  - G_{\alpha,n}));$$
\item for each $\alpha \in \omega_1'$, $n < \omega_0$,
$G_{\alpha+1,n} \cap G_\alpha = G_{\alpha,n}$.
\endroster

We proceed as follows.

First, choose an arbitrary chain
$(G_{0,n} \mid n < \omega_0)$ of finite subgroups with union $G_0$.

For $\alpha \in \omega'$, choose an arbitrary chain
$(\tilde G_{\alpha,m} \mid m < \omega_0)$ of finite subgroups with union
$G_\alpha$.
For each $m < \omega_0$, define $\tilde G_{\alpha+1,m}:=
\langle \tilde G_{\alpha,m}, g_\alpha \rangle.$  Choose any
increasing function $\omega_0 \to \omega_0$, $n \mapsto m_n$, such that the
chain
$(\tilde G_{\alpha+1,m_n}^\bullet \mid n < \omega_0)$ of finite subgroups
of $G_0$
with union $G_0$ is {\it strictly} increasing.  For  $n< \omega_0$, define
$$G_{\alpha+1,n}:= \tilde G_{\alpha+1,m_n}  \text{ and }
G_{\alpha,n}:= G_\alpha \cap G_{\alpha+1,n}.$$
Notice that (4) holds, and that we may assume that (3) holds.  Also
$G_{\alpha,n} \supseteq \tilde G_{\alpha,m_n},$ so
$$G_{\alpha+1,n} = \tilde G_{\alpha+1,m_n}
= \langle \tilde G_{\alpha,m_n}, g_\alpha\rangle
\subseteq \langle G_{\alpha,n},g_\alpha \rangle
\subseteq G_{\alpha+1,n}.$$

For $\alpha < \omega_1$, and $n < \omega_0$, if $G_{\alpha,n}$ is defined, set
$G_{\alpha+1,n}:= \langle G_{\alpha,n},g_\alpha \rangle$; this completes
the recursive definition of the family
$(G_{\alpha,n} \mid \alpha < \omega_1, n < \omega_0)$, and we see that (1)
and (2)
also hold.

For $\alpha < \omega_1$, let $Y_{\alpha}$ denote the $G_{\alpha}$-tree of
Definition~4.4 corresponding to the chain \linebreak $(G_{\alpha,n} \mid n <
\omega_0)$.

We define the {\it plane}, denoted $\Bbb R^2$, to be
the two-dimensional CW-complex with vertices
$(v_{m,n} \mid (m,n) \in \Bbb Z^2)$, and edges
$(x_{m,n}, y_{m,n} \mid  (m,n) \in \Bbb Z^2)$, and faces  \linebreak
$(f_{m,n} \mid  (m,n) \in \Bbb Z^2)$, and incidence relations given by, for
$(m,n) \in \Bbb Z^2$,
$$\iota x_{m,n} = \tau x_{m-1,n} = \iota y_{m,n} = \tau y_{m,n-1} = v_{m,n}$$
and $f_{m,n}$ is attached along the path
$x_{m,n}, y_{m+1,n}, x_{m,n+1}^{-1}, y_{m,n}^{-1}$.

We define the {\it semi-infinite strip}, denoted $[0,1]  \times  \Bbb R^+$,
to be
the subcomplex of $\Bbb R^2$ with vertices
$(v_{0,n}, v_{1,n} \mid n < \omega_0)$, edges
$(x_{0,n}, y_{0,n}, y_{1,n} \mid n < \omega_0)$, and faces
$(f_{0,n} \mid  n < \omega_0)$.

Notice there are two distinguished subcomplexes of
$[0,1]  \times \Bbb R^+$ which are isomorphic to $\Bbb R^+$,
and will be denoted $\{0\} \times \Bbb R^+$ and $\{1\}\times \Bbb R^+$.

Let $\alpha < \omega_1$.  We construct a two-dimensional $G_{\alpha+1}$-space
$M_{\alpha}$, for which \linebreak $[0,1] \times \Bbb R^+$ is a
$G_{\alpha+1}$-trans\-ver\-sal, and the
$G_{\alpha+1}$-stabilizer of $v_{0,n}$, $x_{0,n}$, $y_{0,n}$, and $f_{0,n}$ is
$G_{\alpha,n}$, while the
$G_{\alpha+1}$-stabilizer of  $v_{1,n}$ and $y_{1,n}$ is $G_{\alpha+1,n}$.  This
completely specifies $M_{\alpha}$.  Notice that $M_{\alpha}$ is the mapping
cylinder
of a $G_{\alpha+1}$-map
$$G_{\alpha+1} \times_{G_\alpha} Y_{\alpha}  \to Y_{\alpha+1}.$$
If $\alpha \in \omega_1'$, then  $Y_{\alpha}$ is a $G_{\alpha}$-subtree of
$Y_{\alpha+1}$.

To construct $X$ we glue together
$$(G\times_{G_{\alpha}} M_{\alpha} \mid \alpha < \omega_1)$$ amalgamating
$$(G\times_{G_{\alpha}} Y_{\alpha}\mid \alpha < \omega_1),$$ as in the
proof of Theorem~2.4.   The image in $X_{\alpha+1} \subseteq X$ of
$f_{0,n}\in M_{\alpha}$ will be denoted $f_{\alpha,n}$.

We denote the set of vertices (resp. edges, resp. faces) of $X$ by $VX$ (resp.
$EX$, resp. $FX$).

Let us recall the notion of a club (= {\bf cl}osed {\bf u}n{\bf b}ounded
subset).
Thus, an {\it $\omega_1$-club} is any subset $S$ of $\omega_1$ such that the
set of the least upper bounds of the nonempty subsets of $S$ is precisely
$S \cup \{\omega_1\}$.   If $S$ is an $\omega_1$-club, then so is
$S\cap\omega_1'$; recall that $\omega_1'$ denotes the set of limit ordinals in
$\omega_1$.

For $\phi \in A_G[[X]]$ and $\alpha \le \omega_1$, if
$\supp(\phi) \cap GX_\alpha \subseteq X_\alpha$,  we say that $\phi$ {\it
respects}
$\alpha$.   It is straightforward to show that the set of ordinals in $\omega_1$
respected by $\phi$ is an $\omega_1$-club.
\qed\enddefinition

\proclaim{5.2 Lemma {\rm (Holt~\cite{\HoltOneEnd})}}  Let $G$ be a periodic
group,
and let $H$, $K$ be proper subgroups of $G$ which generate $G$.  Let $X$ be the
$G$-graph with vertex set the disjoint union of $G/H $ and $G/K$, and edge set
$G/(H\cap K)$, with $g(H\cap K)$ joining $gH$ to $gK$, for each $g\in G$.
Then $X$ is
connected, and deleting the two vertices $H$ and $K$, and the one edge $H
\cap K$,
leaves a connected space.
\endproclaim

\demo{Proof} Collapsing all the edges of $X$ leaves a transitive $G$-set, and
one of the points is fixed by $H$ and $K$.  Since $H$ and $K$ generate $G$, this
point is fixed by $G$, so forms a $G$-orbit.  Thus we have only one point, so
$X$ is connected.

Let $Y$ be the subgraph of $X$ obtained by deleting the vertices $H$ and
$K$ and all
their incident edges.  It suffices to show that $Y$ is a connected graph.
Since $X$
is connected, it suffices to show that each $X$-neighbour of $H$  is
$Y$-connected to
each $X$-neighbour of $K$.   Thus, let $h \in H - K$, and $k \in K- H$;  it
suffices
to show that the vertices $hK$ and $kH$ are $Y$-connected.

Let $L = \langle hk^{-1} \rangle =  \langle kh^{-1} \rangle$.  Since $G$ is
periodic,
$L$ is finite.

We consider the action of $L$ on $X$.  Let $m$ and $n$ denote the orders of the
$L$-orbits of the vertices $H$ and $K$, respectively.  By symmetry, we may
assume
that $m \ge n$.  Notice that  $L$  is not contained in $K$, so $n \ge 2$.

Let $g = kh^{-1}$.  In $X$, there is an edge $H\cap K$ joining $H$ to $K$, and
an edge  $k(H\cap K)$ joining $kH = gH$ to $kK = K$.  Applying powers of
$g$ to these,
we get a path in $X$ with vertices
$$H, K, gH, gK, g^2H, g^2K, \ldots, g^{n-1}H, g^{n-1}K.$$
By the definition of $m$ and $n$, these $2n$ vertices are all distinct, so,
on deleting the first two, we get a path in $Y$ connecting $gH = kH$ to
$g^{n-1}K = g^{-1}K = hK$.
\qed \enddemo

\proclaim{5.3 Theorem {\rm (Holt~\cite{\HoltOneEnd})}}  If $G$ is locally
finite,
and $\vert G \vert = \aleph_1$, then $\H^1(G,AG) = 0$.
\endproclaim

\demo{Proof}  Let $X$ be as in Definitions~5.1.

Consider any $\phi \in \Z^1(A_G[[X]])$. Thus $\supp(\phi)$ is a collection
of edges
of $X$, with only finitely many in each $G$-orbit.  A subset of $X$ which meets
(that is, has nonempty intersection with) $\supp(\phi)$ is said to be {\it
broken by}
$\phi$.   Since $\phi$ is a 1-cocycle, we get 0 if we sum, in $A$, the
$\phi$-labels,
with the appropriate signs, around any face, or along any closed path in the
$1$-skeleton, since $X$ is simply-connected.  Thus there is a well-defined {\it
$\phi$-sum} from any vertex to any other vertex.

Consider any $\alpha \in \omega_1'$ such that $\phi$ respects $\alpha$ as in the
last paragraph of Definitions~5.1.

{From} Definitions~5.1, there is a cellular
$G_{\alpha+1}$-map $M_{\alpha+1} \to X$, so $\phi$ induces an element
$\phi_{\alpha+1} \in A[[M_{\alpha+1}]]$.  Since the $G$-stabilizers for
$X$ are finite, $\phi_{\alpha+1}$ lies in
$A_{G_{\alpha+1}}[[M_{\alpha+1}]]$,   Moreover, $\phi_{\alpha+1}$ respects
$\alpha$ in the obvious sense, since
$Y_{\alpha}$ is mapped to $X_\alpha$, by construction.

There exists $n_0 < \omega_0$ such that
$\supp(\phi_{\alpha+1}) \cap G_{\alpha+1}x_{0,0} \subseteq
G_{\alpha+1,n_0}x_{0,0}$.
This means that, for $g \in G_{\alpha+1}$, if $\phi_{\alpha+1}$ breaks
$gx_{0,0}$, then the terminal vertex of $gx_{0,0}$ lies in
$G_{\alpha+1,n_0}v_{1,0}$.

Now consider any $n$  such that $n_0 < n < \omega_0$.

Let $p_{n}$ denote the reduced open path in the tree $\{0\} \times \Bbb
R^+$ from
$v_{0,0}$ to $v_{0,n}$, and  $e \cdot p_{n}$ the open path obtained by
concatenating $e\colon = x^{-1}_{0,0}$ and $p_{n}$.  We are interested in the
$G_{\alpha+1,n}$-graph $Z$ generated by the closure
$\overline{e\cdot p_{n}}$.

Consider any $g \in G_{\alpha+1,n}$.

If $\phi_{\alpha+1}$ breaks $gp_{n}$ then $g\in G_\alpha$, since
$\phi_{\alpha+1}$
respects $\alpha$, and hence  $g \in G_{\alpha,n}$, so
$gv_{0,n} = v_{0,n}$.  That is, if $\phi$ breaks $gp_{n}$, then the terminal
vertex of $gp_{n}$ is $v_{0,n}$.

We apply Lemma~5.2 to the graph $Y'$ obtained by taking  $H^- =
G_{\alpha+1,0}$ and
$K^- = G_{\alpha,n_0}$, so
$\langle H^-, K^- \rangle = \langle G_{\alpha,n_0}, g_\alpha\rangle =
G_{\alpha+1,n_0}$.  By Lemma~5.2, deleting $g^{-1}K^-$ from $Y'$
leaves a connected space containing $G_{\alpha+1,n_0}/H^-$, where we include
the trivial case where $g^{-1}K^-$ does not belong to $Y'$.

Let $Y$ be the $G_{\alpha+1,n_0}$-subspace of $M_{\alpha+1}$ generated by
$\overline{e\cdot p_{n}}$, and consider the
$G_{\alpha+1,n_0}$-map from $Y$ to $Y'$ which assigns $v_{1,0}$ to
$H^-$, $v_{0,n}$ to $K^-$, and $e\cdot p_{n}$ to $H^- \cap K^-$.  It
follows that deleting $g^{-1}K^-\overline{p_{n}}$ from $Y$ leaves a
connected space
containing $G_{\alpha+1,n_0}v_{1,0}$.

Suppose $g \not\in G_{\alpha+1,n_0}$, so $g \in G_{\alpha+1,n} -
G_{\alpha+1,n_0}$.
Then
$\supp(\phi_{\alpha+1}) \cap gY \subseteq K^-p_{n}$.  Hence, on deleting
$\supp(\phi_{\alpha+1})$ from the $1$-skeleton of
$M_{\alpha+1}$,  one of the resulting components contains
$gG_{\alpha+1,n_0}v_{1,0}$.

Next, we apply Lemma~5.2 to the graph $Z'$ obtained by taking $H =
G_{\alpha+1,n_0}$
and   $K = G_{\alpha,n}$, so $\langle H, K \rangle = G_{\alpha+1,n}$.  We
conclude
that
$$Z' - (\{H\} \cup \{K\} \cup  (H\cup K)/(H\cap K))$$
is a connected graph.

Let $Z$ be the $G_{\alpha+1,n}$-subspace of $M_{\alpha+1}$ generated by
$\overline{e\cdot p_{n}}$, and consider the map of $G_{\alpha+1,n}$-spaces
from $Z$
to $Z'$ which assigns $v_{1,0}$ to $H$,  $v_{0,n}$ to $K$, and $e\cdot
p_{n}$ to $H
\cap K$.  There is induced a surjective map
$$Z - (\{v_{1,0}\} \cup \{v_{0,n}\} \cup (H\cup K)(e\cdot p_{n}))
\to  Z' - (\{H\} \cup \{K\} \cup  (H\cup K)/(H\cap K)).$$

Notice that $\phi_{\alpha+1}$ breaks only edges of $Z$ which lie in $He
\cup Kp_{n}$,
so there is a map from
$Z - (\{v_{1,0}\} \cup \{v_{1,n}\} \cup (H\cup K)(e\cdot p_{n}))$ to the set of
components of the $1$-skeleton of $M_{\alpha+1} -\supp(\phi_{\alpha+1})$.
Moreover,
we have  seen that each subset $gHv_{1,0}$ maps to a component of the
$1$-skeleton of
$M_{\alpha+1} -\supp(\phi_{\alpha+1})$.  Thus the map factors through
$Z' - (\{H\} \cup \{K\} \cup  (H\cup K)/(H\cap K))$, which is connected, so
maps to a single component.  Hence some component $X'$ of the $1$-skeleton of
$M_{\alpha+1} -\supp(\phi_{\alpha+1})$ contains
$(\langle H,K \rangle - H)v_{1,0}$, that is,
$(G_{\alpha+1,n} - G_{\alpha+1,n_0})v_{1,0}$.

Since $n > n_0$ was arbitrary, all of $(G_{\alpha+1} -
G_{\alpha+1,n_0})v_{1,0}$
is contained in $X'$.  Thus, for any path between any two elements of
$(G_{\alpha+1} - G_{\alpha+1,n_0})v_{1,0}$, the $\phi_{\alpha+1}$-sum, and
hence the
$\phi$-sum, is 0.

Let $\psi_{\alpha+1} \in \C^0(A[[X_{\alpha+1}]])$ be defined on each vertex
$v$ as the $\phi$-sum along any path from any vertex of
$(G_{\alpha+1} - G_{\alpha+1,n_0})v_{1,0}$ to $v$. Then
$\psi_{\alpha+1} \in \C^0(A_{G_{\alpha+1}}[[X_{\alpha+1}]])$ and
$\phi\vert_{X_{\alpha+1}} = \partial^*(\psi_{\alpha+1})$.

The $\alpha<\omega_1$ which are respected by $\phi$ converge to $\omega_1$,
and it
follows that the corresponding (unique) $\psi_{\alpha+1}$ converge to an
element
$\psi \in \C^0(A_G[[X]])$ such that $\phi = \partial^*(\psi)$, so
$\phi \in \B^1(A_G[[X]])$.

Hence $\H^1(A_G[[X]]) = 0$, so $\H^1(G,AG) = 0$.
\qed \enddemo

Up until now, in this section, we have given a straightforward
topological translation of~\cite{\HoltOneEnd}, which we felt illuminated the
arguments.   We now come to a new result.

\proclaim
{5.4 Theorem} If $G$ is locally finite, and $\vert G \vert = \aleph_1$, then
$\H^2(G, AG)$ has a subgroup isomorphic to $A^{\aleph_1}$, and
$\vert\H^2(G, AG)\vert=\vert A^G\vert=\vert A\vert^{\aleph_1}$.
\endproclaim

\demo{Proof}  Let $X$ be as in Definitions~5.1.

Let $\alpha \in \omega_1'$.

Let $M'_\alpha$ denote the $G_\alpha$-space with $G_\alpha$-transversal
$[-1,0] \times \Bbb R^+$,  where the $G_\alpha$-sta\-bi\-liz\-ers of
$v_{-1,n}$, $y_{-1,n}$, $x_{-1,n}$ and $f_{-1,n}$ are
$G_{\alpha,n}^\bullet (= G_0 \cap G_{\alpha,n})$, while the
$G_{\alpha}$-sta\-bi\-liz\-ers of  $v_{0,n}$ and $y_{0,n}$ are
$G_{\alpha,n}$.  This
completely specifies $M'_{\alpha}$.

Notice that $M'_{\alpha}$ is the mapping
cylinder of an injective $G_{\alpha}$-map
$G_{\alpha} \times_{G_0} Y_{\alpha}^\bullet  \to Y_{\alpha},$
where $Y_{\alpha}^\bullet$ is the $G_{0}$-tree of Definition~4.4
corresponding to
the chain $(G_{\alpha,n}^\bullet \mid n < \omega_0)$, a $G_{0}$-subtree
of $Y_{\alpha}$.

We have a cellular $G_\alpha$-map $Y_\alpha \to X_\alpha$, and, similarly,
by Corollary~2.3, we can construct a cellular $G_0$-map $Y_{\alpha}^\bullet
\to Y_0$
between spaces of type $\uE G_0$.  These two maps can be extended to a
$G_\alpha$-map
$M'_\alpha \to X_\alpha$.

Notice that $Y_\alpha$ is contained in both $M_{\alpha}$ and
$M'_{\alpha}$, and the map $Y_\alpha \to X_\alpha$ has been extended to
$M_{\alpha} \to X_{\alpha+1}$ and to $M'_\alpha \to X_\alpha$.

For each $n < \omega_0$, let $\lambda_{\alpha,n}'$ denote the least ordinal such
that the $X_{\lambda_{\alpha,n}'+1}$ contains the image of the face
$f_{-1,n}$, under
the map $M'_\alpha \to X_\alpha$.  Notice that $\lambda_{\alpha,n}'+1 < \alpha$,
since
$\alpha$ is a limit ordinal. Choose a strictly ascending sequence
$(\lambda_{\alpha,n} < \alpha \mid n < \omega_0)$ with limit $\alpha$,
such that
$\lambda_{\alpha,n} > \max\{\lambda_{\alpha,i}' \mid 0 \le i \le n\}$.
Let $$h_{\alpha,n}:=  g_{\lambda_{\alpha,n}} \in
G_{\lambda_{\alpha,n}+1} - G_{\lambda_{\alpha,n}}.$$

Let $Q =
\{h_{\alpha,n} f_{\alpha,n} \mid
(\alpha,n) \in \omega_1'\times \Bbb N\}$,  so
$$A[[Q]] \subseteq A_G[[FX]] = \C^2(A_G[[X]]) = \Z^2(A_G[[X]]).$$

Recall that a subset $S$ of $\omega_1$ is said to be {\it
$\omega_1$-stationary}
if $S$ meets each $\omega_1$-club.  Let $\Cal P$ be a partition of
$\omega_1'$ into $\aleph_1$ subsets, each being $\omega_1$-stationary;
see, for example,~\cite{\Jech,~Lemma~7.6,~p.59}.

Consider any  $\phi \in A^{\Cal P}$, $[\alpha] \mapsto \phi[\alpha]$.

Define
$$\phi^\dag\colon= \sum\limits_{(\alpha,n) \in \omega_1'\times \Bbb N}
\phi[\alpha].h_{\alpha,n}f_{\alpha,n} \in A[[Q]]
\subseteq
\Z^2(A_G[[X]]).$$
It is easy to see that $\phi^\dag$ respects all the ordinals in $\omega_1$.

We claim that if $\phi^\dag \in \B^2(A_G[[X]])$, then $\phi = 0$.

Suppose not, so there exists $\psi \in \C^1(A_G[[X]]) = A_G[[EX]]$, such that
$\partial^*\psi = \phi^\dag$, and there exists $p \in \Cal P$ such that
$\phi(p) \ne 0$. We shall obtain a contradiction.

Since $Y_0$ is countable, there exists  $\mu < \omega_1$ such that
$\supp(\psi) \cap GY_0 \subseteq G_\mu Y_0$.

Since $p$ is $\omega_1$-stationary, and the set of ordinals respected by
$\psi$ is an
$\omega_1$-club, there exists $\alpha \in p$ such that $\alpha > \mu$ and $\psi$
respects $\alpha$.  Notice $\phi[\alpha] = \phi(p) \ne 0$.

Now $\phi^\dag$ and $\psi$ induce elements $\phi_{\alpha+1}$ and
$\psi_{\alpha+1}$, respectively, in $A_{G_{\alpha+1}}[[M_{\alpha}]]$, and,
moreover,  $\partial^*\psi_{\alpha+1} =  \phi_{\alpha+1}$ in
$A_{G_{\alpha+1}}[[M_{\alpha}]]$.  Notice
$\phi_{\alpha+1} =   \sum\limits_{n < \omega_0} \phi(p).h_{\alpha,n}f_{0,n}$.

Let $\xi:=  - \sum\limits_{n < \omega_0} \phi(p).h_{\alpha,n}y_{0,n} \in
A_{G_{\alpha}}[[EY_{\alpha}]] \subseteq A_{G_{\alpha+1}}[[M_{\alpha}]]$.  Then
$\partial^*(\xi) =  \phi_{\alpha+1}$ in $A_{G_{\alpha+1}}[[M_{\alpha}]]$.  Thus
$\partial^*(\psi_{\alpha+1} - \xi) = 0$.  Moreover, $\psi_{\alpha+1} - \xi$
respects $\alpha$, since both $\xi$ and  $\psi_{\alpha+1}$ respect
$\alpha$.  By the
proof of Theorem~5.3,  there exists $\nu < \alpha$
such that the $(\psi_{\alpha+1} - \xi)$-sum along any path
in $Y_{\alpha}$ between any two vertices in
$(G_{\alpha} - G_\nu)v_{0,0}$ is zero.

Now $\phi^\dag$ and $\psi$ induce elements $\phi_{\alpha}$ and
$\psi_{\alpha}$ in
$A_{G_{\alpha}}[[M'_{\alpha}]]$; moreover,
$\partial^*\psi_{\alpha} =  \phi_{\alpha}$ in
$A_{G_{\alpha}}[[M'_{\alpha}]]$, and
$\psi_{\alpha+1}$ and $\psi_\alpha$ agree on $Y_\alpha$.
Also, $\xi$ can be viewed as an element of $A_{G_{\alpha}}[[M'_{\alpha}]]$,
so the
$(\psi_{\alpha} - \xi)$-sum along any path in $Y_{\alpha}$ between any two
vertices in  $(G_{\alpha} - G_\nu)v_{0,0}$ is zero.

There exists $\kappa < \alpha$ such that $\psi_\alpha$ vanishes on
$(G_{\alpha} - G_\kappa)x_{-1,0}$.

Choose  $n < \omega_0$ such that $\lambda_{\alpha,n}$ is greater than
$\mu$, $\nu$ and $\kappa$.  Choose \linebreak
$g \in G_{\alpha,n+1}^\bullet - G_{\alpha,n}^\bullet$.
Thus $g \in G_0$, $g$ fixes $x_{-1,n+1}$, and $g$ moves $y_{0,n}$.

In $\Bbb Z[M_\alpha]$,
$$\partial(\sum\limits_{i=0}^n f_{-1,i})
=  x_{-1,0} - x_{-1,n+1} + \sum\limits_{i=0}^n (y_{0,i}-y_{-1,i}).$$
Hence $$\partial(h_{\alpha,n}(1-g)\sum\limits_{i=0}^n f_{-1,i})
= h_{\alpha,n}(1-g)(x_{-1,0}+\sum\limits_{i=0}^n (y_{0,i}-y_{-1,i})).$$
We can view $\psi_{\alpha}$ as an additive map $\Bbb Z[EM_\alpha]\to A$,
and apply
it to the foregoing equation, to get
$$\align
\phi_{\alpha}(h_{\alpha,n}(1-g)&\sum\limits_{i=0}^n f_{-1,i}) \\ =
&\psi_{\alpha}(h_{\alpha,n}(1-g)(x_{-1,0}+\sum\limits_{i=0}^n(y_{0,i}-y_{-1,
i}))).
\endalign$$

Notice that $h_{\alpha,n}G_{\lambda_{\alpha,n}} \cap G_{\lambda_{\alpha,n}} =
\emptyset$.

Since $h_{\alpha,n}, h_{\alpha,n}g \not\in G_{\lambda_{\alpha,i}'}$,
for $0 \le i \le n$,
we see $\phi_{\alpha}(h_{\alpha,n}(1-g)\sum\limits_{i=0}^n f_{-1,i}) = 0$.

Also, $h_{\alpha,n}, h_{\alpha,n}g \notin G_{\mu}$, so
$\psi_{\alpha}(h_{\alpha,n}(1-g)y_{-1,i})=0$.

Also, $h_{\alpha,n}, h_{\alpha,n}g \notin G_{\kappa}$, so
$\psi_{\alpha}(h_{\alpha,n}(1-g)x_{-1,0}) = 0$.

It follows that
$\psi_{\alpha}(h_{\alpha,n}(1-g) \sum\limits_{i=0}^n y_{0,i}) = 0.$

But $h_{\alpha,n}, h_{\alpha,n}g \notin G_{\nu}$, so
$(\psi_{\alpha}-\xi)(h_{\alpha,n}(1-g)\sum\limits_{i=0}^n y_{0,i}) =0.$
Hence  $$\xi(\sum\limits_{i=0}^n h_{\alpha,n}(1-g)y_{0,i}) = 0.$$  Now
$\xi(h_{\alpha,n}y_{0,n}) = -\phi(p)$, and $\xi$ vanishes on all other
summands  because \linebreak $h_{\alpha,n}G_0 \cap G_{\lambda_{\alpha,i}+1} =
\emptyset$ for
$0 \le i \le n-1$, and $gy_{0,n} \ne y_{0,n}$.  Hence $\phi(p) = 0$, which is a
contradiction.

Since $\vert \Cal P \vert = \aleph_1$,  we have an embedding of $A^{\aleph_1}$
in $\H^2(G,AG)$.
\qed
\enddemo

\head 6. Cohomology of directed unions \endhead

In this section, we recall some known results about cohomology for well-ordered
directed unions, with special emphasis on abelian groups.

\definition
{6.1 Notation}
We let $(\P(G),\partial)$ denote the bar resolution for $G$, and let $\P_n(G)$
denote its $n$th component, for each $n \in \Bbb Z$.  Thus $(\P(G),\partial)$
is a free $\Bbb ZG$-resolution of $\Bbb Z$, and, for $n \ge 0$, $\P_n(G)$ has as
$\Bbb Z$-basis the Cartesian power $G^{n+1}$, with $G$ acting by left
multiplication
on the first coordinate, and, for $n \ge 1$,
$$\align
\partial_n(g_0,&\ldots,g_n):=\\
&\sum_{i=0}^{n-1} (-1)^i(g_0,\ldots, g_{i-1}, g_ig_{i+1},g_{i+2},\ldots, g_n)
+ (-1)^{n}(g_0,\ldots,g_{n-1}).\endalign $$
As usual, if $n \le -1$,  $\P_n(G) = 0$,  and, if $n \le 0$, $\partial_{n}
= 0$.
\qed\enddefinition

The following is a degenerate case of the cohomology spectral sequence for
well-ordered directed unions; see, for example~\cite{\Robinson,~Section~3}.

\proclaim
{6.2 Lemma {\rom{(Robinson \cite{\Robinson, Proposition~1})}}}
Let $n \in \Bbb N$,  let $M$ be a $\Bbb ZG$-mod\-ule, let $\beta$ be a
limit ordinal, and let $(G_\alpha \mid \alpha \le \beta)$ be a continuous
chain of
subgroups of $G$.  If $\H^{n-1}(G_\alpha,M)= 0$ for all $\alpha < \beta$, then
$\H^n(G_\beta,M) = \lim\limits_{\overleftarrow {\alpha < \beta}}
\H^n(G_\alpha,M)$.
\endproclaim

\demo{Proof}  For each $\alpha \le \beta$, view the bar resolution
$\P(G_\alpha)$
as a $\Bbb Z$-subcomplex of $\P(G_\beta)$, so
$(\P_n(G_\alpha)\mid \alpha \le \beta)$ is a continuous chain.

We want to show that the natural map
$$\H^n(G_\beta,M) \to \lim\limits_{\overleftarrow {\alpha < \beta}}
\H^n(G_\alpha,M)\tag"(6.3)"$$
is bijective.

We begin by showing it is injective.

Consider any element $\xi$ of the kernel of (6.3).  Then $\xi$ is represented
by an $n$-cocycle  $\phi_\beta\colon \P_n(G_\beta) \to M$, so $\phi_\beta$ is
$\Bbb ZG_\beta$-linear, and $\phi_\beta \circ \partial_{n+1} = 0$.

Consider any $\alpha \le \beta$.  Let $\phi_{\beta,\alpha}$ denote the
restriction
of $\phi_\beta$ to $\P_n(G_\alpha)$.  We shall construct, transfinitely, a
continuous
directed system of maps
$$(\psi_{\beta,\alpha}\colon \P_{n-1}(G_\alpha) \to M \mid \alpha \le
\beta)$$ such
that $\psi_{\beta,\alpha}$ is $\Bbb ZG_\alpha$-linear, and
$\psi_{\beta,\alpha} \circ \partial_n = \phi_{\beta,\alpha}$.
It will then follow that \linebreak
$\phi_\beta = \phi_{\beta,\beta} = \psi_{\beta,\beta} \circ \partial_n$ is a
coboundary, that $\xi = 0$, and that (6.3) is injective.

If $\alpha < \beta$, then $\phi_{\beta,\alpha}$ represents an element
$\xi_\alpha$
of  $\H^n(G_\alpha,M)$, and, since $\xi$ lies in the kernel of (6.3),
$\xi_\alpha =
0$, so there exists a $\Bbb ZG_\alpha$-linear map
$$\psi_\alpha\colon \P_{n-1}(G_\alpha) \to M$$ such that
$\psi_\alpha \circ \partial_n = \phi_{\beta,\alpha}$.

Let $\psi_{\beta,0} = \psi_0$.

If $\alpha < \beta$ and we have constructed $\psi_{\beta,\alpha}$, then we
construct $\psi_{\beta,\alpha+1}$ as follows.  Let $\psi_{\alpha+1,\alpha}$
denote
the restriction of $\psi_{\alpha+1}$ to $P_{n-1}(G_{\alpha})$.
Since $\psi_{\alpha+1} \circ \partial_n =  \phi_{\beta,\alpha+1}$, we see that
$\psi_{\alpha+1,\alpha} \circ \partial_n =  \phi_{\beta,\alpha}
= \psi_{\beta,\alpha} \circ  \partial_n$. Thus
$\psi_{\beta,\alpha} - \psi_{\alpha+1,\alpha}$ is an $(n-1)$-cocycle, so
represents an
element of $\H^{n-1}(G_\alpha, M)$.  By hypothesis, this element is zero,
so there
exists a $\Bbb ZG_\alpha$-linear map $\mu_\alpha\colon \P_{n-2}(G_\alpha)
\to M$ such
that $\psi_{\beta,\alpha} - \psi_{\alpha+1,\alpha}= \mu_\alpha\circ
\partial_{n-1}$.
Let $\mu_{\alpha,\alpha+1}\colon \P_{n-2}(G_{\alpha+1}) \to M$ denote the
unique $\Bbb ZG_{\alpha+1}$-linear map which is $\mu_{\alpha}$ on
$\P_{n-2}(G_\alpha)$, and is zero on
$G_{\alpha+1}^{n-1} - (G_{\alpha+1} \times G_{\alpha}^{n-2})$.
Now define
$$\psi_{\beta,\alpha+1}\colon
= \psi_{\alpha+1}  + \mu_{\alpha,\alpha+1}\circ \partial_{n-1}.$$
By construction, $\psi_{\beta,\alpha+1}$ acts on $\P_{n-1}(G_\alpha)$ as
$\psi_{\alpha+1,\alpha}  + \mu_{\alpha}\circ \partial_{n-1} =
\psi_{\beta,\alpha}$.
Also $\psi_{\beta,\alpha+1}\circ \partial_n = \psi_{\alpha+1}\partial_n + 0 =
\phi_{\beta,\alpha+1}$.

If $\alpha$ is a limit ordinal in $\beta+1$ and we have constructed
$(\psi_{\beta,\alpha'} \mid \alpha' < \alpha)$, then the latter has a
direct limit
which we take to be $\psi_{\beta,\alpha}$.

This completes the construction, so (6.3) is injective.

A similar, but easier,  argument shows that (6.3) is surjective, and here the
hypothesis on $\H^{n-1}$ is not needed.
\qed\enddemo

A special case of this result appeared in the penultimate paragraph of the
proof of
Theorem~5.3.  We can now say even more.

\proclaim{6.4 Theorem {\rm (Holt~\cite{\HoltOneEnd})}}  If $G$ is locally
finite,
and $\arank(G) \ne 0$, then \linebreak $\H^1(G,AG) = 0$.
\endproclaim

\demo{Proof} By (1.3), we may assume that $\arank(G) \ge 1$, and we proceed by
induction on $\arank(G)$.  If $\arank(G) = 1$, the assertion holds by
Theorem~5.3.
Thus we may assume that $\arank(G) \ge 2$ and the result holds for smaller
groups.
There exists an ordinal $\beta$ and a continuous chain
 $(G_\alpha \mid\alpha\le\beta)$ of subgroups of $G$ with $G_\beta = G$, and
$1 \le \arank(G_\alpha) < \arank(G)$ for all $\alpha < \beta$.

By (1.3) and the induction hypothesis, $\H^0(G_\alpha, AG) = \H^1(G_\alpha,
AG) = 0$
for all $\alpha < \beta$, since $AG$ is an induced $\Bbb ZG_\alpha$-module.  By
Lemma~6.2, $\H^1(G,AG) = 0$.
\qed \enddemo

\definition
{6.5 Remark} In light of Theorem~3.10, the method of proof of Theorem~6.4
shows that
$$\text{``If } n \in \Bbb N, \,G \text{ is locally finite, and }
n \ne \arank(G)+1, \text{ then } \H^n(G,\Bbb ZG\otimes_{\Bbb Z} -) =
0,\text{"}$$ is
equivalent to
$$\text{``If } n \in \Bbb N, \, G \text{ is locally finite, and }
n = \arank(G), \text{ then }\H^n(G,\Bbb ZG\otimes_{\Bbb Z} -) = 0.\text{"}\qed$$
\enddefinition

\definition
{6.6 Definitions} Let $H$ be a subgroup of $G$, and let $L$, $M$ be
 $\Bbb ZG$-modules.

Let $N_G(H)$ denote the normalizer of $H$ in $G$, and let $C_G(H)$ denote the
centralizer of $H$ in $G$.

There is a natural action of $N_G(H)$ on $\Ext_{\Bbb Z H}^*(L,M)$, where
$L$, $M$ are
viewed as $\Bbb ZH$-modules by restriction.  Perhaps the simplest way to
define this action is as follows.  Let $(\P,\partial)$ be a projective
$\Bbb ZG$-resolution of $L$, and view $(\P,\partial)$ as a projective
$\Bbb ZH$-resolution of $L$.  Then any  $g \in N_G(H)$ gives rise to an
action on $\Hom_{\Bbb ZH}(\P,M)$, $\phi \mapsto \phi^g$, where
$\phi^g(p) = g^{-1}(\phi(gp))$.   Since $g$ normalizes $H$, we see that
$\phi^g$ is
$\Bbb ZH$-linear.  The $N_G(H)$-action respects cocycles and coboundaries, so
induces an action on the cohomology, $\Ext_{\Bbb Z H}^*(L,M)$.

In the foregoing, if $\phi$ is $\Bbb ZG$-linear, then $\phi^g = \phi$.  It
follows
that the image of the restriction map
$$\Ext_{\Bbb Z G}^*(L,M) \to \Ext_{\Bbb Z H}^*(L,M)$$
is contained in the set $\Ext_{\Bbb Z H}^*(L,M)^{N_G(H)}$ of points fixed by
$N_G(H)$.

Let $g$ be an element of $C_G(H)$.  Left multiplication by $g$ determines a
$\Bbb ZH$-linear endomorphism of $M$ and we denote it
by $g\cdot\vert_M$.  Similarly, $g\cdot\vert_{\P} \in \End_{\Bbb ZH}(\P)$
is a lift of $g\cdot\vert_{L} \in \End_{\Bbb ZH}(L)$.  Now $\Ext_{\Bbb Z
H}^*(L,M)$ is
a bimodule over $\End_{\Bbb ZH}(M)$ and $\End_{\Bbb ZH}(L)$, and we see that the
action of $g$ on $\Ext_{\Bbb Z H}^*(L,M)$ is given by
$\eta \mapsto (g\cdot\vert_M)^{-1} \circ \eta \circ (g\cdot\vert_L)$.

If we restrict to the case where $L = \Bbb Z$ with trivial $G$-action, we
see that
the image of the restriction map $\H^*(G,M) \to \H^*(H,M)$ is contained in
$\H^*(H,M)^{N_G(H)}$, and hence in $\H^*(H,M)^{C_G(H)}$, and here $C_G(H)$
acts by
left multiplication on $M$. \qed
\enddefinition

The following is now an immediate consequence of Lemma~6.2

\proclaim
{6.7 Corollary}  Let $n \in \Bbb N$,  let $M$ be a $\Bbb ZG$-mod\-ule, let
$\beta$ be
a limit ordinal, and let $(G_\alpha \mid \alpha \le \beta)$ be a continuous
chain of
subgroups of $G$ with $G_\beta = G$.  If, for each $\alpha < \beta$,
$\H^n(G_\alpha,M)^{N_G(G_\alpha)} = 0$ and $\H^{n-1}(G_\alpha,M) = 0$, then
$\H^n(G,M) = 0$.\qed
\endproclaim

We record consequences for abelian groups which seem to be new.

\proclaim
{6.8 Theorem}  Let $G$ be an abelian group,  $\lambda$ an ordinal, and $\Delta$
a $G$-set with stabilizers of $\aleph$-rank strictly less than $\lambda$.
\roster
\item For each $n \in \Bbb N$, if $\arank(G) \ge \lambda+n$  then
$\H^n(G, A\Delta) = 0$.
\item If $\arank(G) \ge \lambda+\omega_0$ then $\H^*(G,A\Delta) = 0$.
\endroster
\endproclaim

\demo{Proof} (1).  We argue by induction on $n$.

If $n = 0$, then all $G$-stabilizers of elements of $\Delta$ have infinite index
in $G$, so $\Delta$ has no finite $G$-orbits.  Thus
$(A\Delta)^G = 0$, that is, $\H^0(G, A\Delta) = 0$.

Now suppose that $n \ge 1$,  and that the result holds for smaller $n$.  Let
$\beta$ denote the least ordinal of cardinality $\rank(G)$, so
$\beta$ is a limit ordinal.  Moreover, there exists a continuous chain of
subgroups
$(G_\alpha \mid \alpha \le \beta)$ such that $G_\beta = G$, and, for each
$\alpha < \beta$,
$$\lambda+n-1 \le \arank(G_\alpha)  < \arank(G).$$

Consider $\alpha < \beta$.

By the induction hypothesis, $\H^{n-1}(G_\alpha,A\Delta) = 0$, so, by
Corollary~6.7,
it remains to show that
$\H^n(G_\alpha,A\Delta)^{N_G(G_\alpha)} =  0$.   Here $N_G(G_\alpha) =
C_G(G_\alpha) = G$, since $G$ is abelian.  Thus, we want to show that
$\H^n(G_\alpha,A\Delta)^G =  0$ where $G$ is acting via multiplication on
$A\Delta$.

Consider any element $\zeta \in \H^n(G_\alpha, A\Delta)$.  We can use the
bar resolution $\P(G_\alpha)$, and represent $\zeta$ by a
$\Bbb ZG_\alpha$-linear map $\phi\colon \P_n(G_\alpha) \to  A\Delta$.  Since
$\vert G_\alpha \vert < \vert G \vert$, we see there is a $G_\alpha$-subset
$\Delta'$
of $\Delta$ such that $A\Delta'$ contains the image of $\phi$, and
$\vert \Delta'\vert < \vert G \vert$.  Let
$\Delta''$ be the complement of $\Delta'$ in $\Delta$.  Then
$A\Delta = A\Delta' \oplus A\Delta''$, so
$$\H^n(G_\alpha, A\Delta)
= \H^n(G_\alpha, A\Delta') \oplus \H^n(G_\alpha, A\Delta''),$$
and, in the corresponding expression $\zeta = (\zeta',\zeta'')$, we have
$\zeta''=0$.

Now $$\vert \{g \in G \mid g\Delta' \cap \Delta' \ne \emptyset\}\vert
< \vert G\vert$$ because the elements of $\Delta$ have $G$-stabilizers of
cardinality
strictly less than $\vert G \vert$.  Hence there exists $g \in G$ such that
$g\Delta' \cap \Delta' = \emptyset$, that is,  $g\Delta' \subseteq \Delta''$, so
$$\H^n(G_\alpha, A\Delta')^g \subseteq \H^n(G_\alpha, A\Delta'').$$
It follows that if $\zeta^g = \zeta$ then $\zeta = 0$.  This proves that
$\H^n(G_\alpha, A\Delta)^G =  0$, and (1) is proved.

(2) follows from (1).
\qed\enddemo

The case where $\Delta = G$ is of interest;  here we can take $\lambda = 0$.

\proclaim
{6.9 Corollary}  Let $n \in \Bbb N$, and let $G$ be an abelian group. Then,
$\H^n(G,AG) = 0$ if $n < \arank(G)+1$. 	\qed
\endproclaim

\proclaim
{6.10 Corollary}  Let $n \in \Bbb N$, and let $G$ be an abelian, torsion group.
\roster
\item {\rm (Holt~\cite{\HoltDimension})}. $\H^n(G,AG) = 0$ if $n\ne
\arank(G)+1$.

\item {\rm (Chen~\cite{\Osofsky,~Corollary~7.6})}.  If $R$ is a nonzero,
$o(G)$-inverting ring, then
\linebreak $\cd_RG = n$ if and only if $\arank(G) = n-1$, that is,
$$\cd_RG = \min\{\arank(G)+1, \infty\}.\qed$$
\endroster
\endproclaim

\demo{Proof}
(1) follows from Theorem~3.10 and Corollary~6.9.
(2) follows from the fact that, if $\cd_R(G) < \infty$, then
$\H^m(G,\Bbb ZG\otimes_{\Bbb Z} -) \ne 0$ for $m = \cd_R(G)$; see the Commentary
on Conjecture~1.7. \qed\enddemo

\head 7. Cardinals, free abelian groups, and $\PHKh\frak{F}$  \endhead

We now recall the hierarchies introduced in~\cite{\Krop}; see~\cite{\KM}
for more
details.

\definition
{7.1 Notation}
Let $\frak{X}$ denote a class of groups.

All the classes of groups that we consider are closed under isomorphism,
for example,
the class $\frak F$ of all finite groups.

We let $\PHKl\frak X$ denote the class of groups whose finitely
generated subgroups all lie in $\frak X$.  For example, if $\frak X$ contains
all finitely generated abelian groups, then $\PHKl\frak X$ contains all abelian
groups.

We let $\PHKh_1\frak{X}$ denote the class of all groups $G$ for which there
exists a finite-di\-men\-sion\-al contractible $G$-complex with all
stabilizers lying
in $\frak X$.   For example, $\PHKh_1\frak{F}$ contains all finitely
generated abelian
groups, since, if $G$ is finitely-generated and abelian, then $G$ has a finite
subgroup $N$ such that $G/N$ is isomorphic to $\Bbb Z^n$ for some $n \in
\Bbb N$, and
thus $G/N$ acts freely on $\Bbb R^n$ preserving a CW-structure.

If $\PHKh_1\frak{X} = \frak{X}$, then $\frak X$ is said to be
{\it $\PHKh_1$-closed}.

We let $\PHKh\frak{X}$ denote the smallest $\PHKh_1$-closed class of groups
which contains $\frak X$.  This class has a hierarchy indexed by the
ordinals, where
for each ordinal $\beta$, we define $\PHKh_\beta\frak{X}$ recursively, by
setting
$$\align
\PHKh_0\frak{X}:&= \frak X,\\
\PHKh_\beta\frak X:&= \PHKh_1\PHKh_{\beta-1} \frak X
\text{ if $\beta$ is a successor ordinal},\\
\PHKh_\beta\frak X:&= \bigcup\limits_{\alpha < \beta} \PHKh_\alpha \frak X
\;\text{
if
$\beta$ is a limit ordinal.} \qed
\endalign$$
\enddefinition

We can use Theorem~2.8 to get new sufficient conditions for membership in
$\PHKh\frak{X}$.

\proclaim
{7.2 Theorem} Let $\frak X$ be a subgroup-closed class of groups, and let
$G \in \PHKl \frak{X}$.
\roster
\item  If $\arank(G)=-1$ then $G \in \PHKh_0\frak X$.
\item If $\arank(G) < \omega_0$ then $G \in \PHKh_1\frak X$.
\item  If $\arank(G) = \omega_0$ then $G \in \PHKh_2\frak X$.
\endroster
\endproclaim

\demo{Proof} (1).  If $\arank(G) = -1$ then $G$ is a finitely generated element
of $\PHKl \frak{X}$, so lies in $\frak{X} = \PHKh_0\frak X$.

(2). If $\arank(G) =  n$ for some $n \in \Bbb N$ then, by Theorem~2.8, $G$
acts on a contractible $(n+1)$-dimensional CW-complex with stabilizers
contained in
finitely generated subgroups of $G$.  Since $G$ lies in  $\PHKl\frak X$, and
$\frak X$ is subgroup closed, we see that these stabilizers lie in $\frak
X$.  Hence
$G \in \PHKh_1\frak X$.

(3).  If $\arank(G) = \omega_0$, then we can write $G$ as the union
of an ascending chain $(G_n \mid n \in \Bbb N)$ of subgroups, such that,
for each $n
\in \Bbb N$,  $\arank(G_n) = n$.  As in Remark~2.5, there exists a
$G$-tree with each cell stabilizer contained in $G_n$, for some $n \in \Bbb
N$, so
lying in $\PHKh_1\frak X$ by (2).  Hence,  $G \in \PHKh_2\frak X$.
\qed
\enddemo

We next consider some necessary conditions for membership in $\PHKh\frak
X$, which
are natural generalizations of~\cite{\KropTwo,~Lemma~1}.

\proclaim
{7.3 Lemma} Let $n \in \Bbb N$, let $R$ be a ring, let
$$0 \to M_n \to M_{n-1} \to \cdots \to M_1 \to M_0 \to M_{-1} \to 0$$
be an exact sequence of $R$-modules, and let $L$ be an $R$-module.

Suppose that $\Ext_R^i(L,M_i) = 0$ for $i = 0, \ldots, n$.  Then
$\Ext^0_R(L,M_{-1}) = 0$.
\endproclaim

\demo{Proof}  Clearly the result holds for $n=0$.   Thus we may assume that
$n \ge 1$, and that the result holds with  $n-1$ in place of $n$.

Let $M'_{n-1}$ denote the cokernel of the map $M_n \to M_{n-1}$, so we have
exact
sequences
$$\align
0 \to M_n \to M_{n-1} \to &M'_{n-1} \to 0\tag"(7.4)"\\
0 \to &M'_{n-1} \to M_{n-2} \to \cdots  \to M_0 \to M_{-1} \to 0\tag"(7.5)".
\endalign$$
Now (7.4) gives rise to a long exact sequence which contains the segment
$$\Ext_R^{n-1}(L, M_{n-1}) \to \Ext_R^{n-1}(L, M'_{n-1}) \to \Ext_R^n(L,M_n).$$
Here the outer terms are zero, by hypothesis,  so the inner term is zero.  The
induction hypothesis can now be applied to (7.5), and we see that
$\Ext^0_R(L,M_{-1}) = 0$.

The result follows.
\qed
\enddemo

We record the contrapositive of the case where $R = \Bbb ZG$, and
$M_{-1} = L = \Bbb Z$ with trivial $G$-action.

\proclaim
{7.6 Corollary} If $n \in \Bbb N$,  and
$$0 \to M_n \to M_{n-1} \to \cdots \to M_1 \to M_0 \to \Bbb Z \to 0$$
is an exact sequence of $\Bbb ZG$-modules, then there exists $i$ such that
$0 \le i \le n$ and $\H^i(G,M_i) \ne 0$.  \qed
\endproclaim

\proclaim
{7.7 Proposition} If $\frak X$ is a class of groups, and $\H^*(G,\Bbb
Z\Delta) = 0$
for every $G$-set $\Delta$ for which all stabilizers lie in $\frak X$, then
$G \not\in \PHKh_{1}\frak X$.
\endproclaim

\demo{Proof}  Suppose that $G \in \PHKh_{1} \frak X$, so there exists a
finite-dimensional, contractible CW-complex $X$ on which $G$ acts with all
stabilizers lying in $\frak X$.  Let $n$ denote the dimension of $X$.  The
augmented
cellular chain complex of $X$,
$$0 \to \C_n(X) \to \C_{n-1}(X) \to \cdots \to\C_1(X) \to \C_0(X) \to \Bbb
Z \to 0,$$
is an exact sequence of $\Bbb ZG$-modules, so, by Corollary~7.6, there exists
$i$ such that $0 \le i \le n$ and $\H^i(G,C_i(X)) \ne 0$.  Thus
$\H^\ast(G,C_i(X)) \ne 0$.  But we can write $C_i(X) = \Bbb Z\Delta$, where
$\Delta$ is the set of $i$-dimensional open cells of $X$, so is a $G$-set with
all stabilizers in $\frak X$.  This contradicts the hypothesis.  \qed
\enddemo

Combining Theorem~6.8(2) and Proposition~7.7, we get the following.

\proclaim
{7.8 Corollary} Let $\frak X$ be a class of groups and $\lambda$ be an ordinal.
Suppose that $G$ is an abelian group in $\PHKh_1\frak X$ such that every
subgroup
of $G$ which lies in $\frak X$ has $\aleph$-rank strictly less than
$\lambda$.  Then $\arank(G) < \lambda+\omega_0$. \qed
\endproclaim

\proclaim
{7.9 Theorem} If $\beta$ is any ordinal, then every abelian group in
$\PHKh_{\beta} \frak F$ has $\aleph$-rank strictly less than $\omega_0\beta$.
\endproclaim

\demo {Proof} We argue by induction on $\beta$.

The result holds for $\beta = 0$ by definition of $\frak F$.

Thus we may assume that $\beta > 0$, and that the result holds for smaller ordinals.

Consider the case where $\beta$ is a limit ordinal.  Here, each abelian group in
$\PHKh_\beta \frak F$ lies in $\PHKh_\alpha \frak F$ for some $\alpha <
\beta$, so,
by the induction hypothesis, is of cardinality strictly less than
$\aleph_{\omega_0\alpha} \le \aleph_{\omega_0\beta}$.

Now consider the case where $\beta$ is a successor ordinal, and write
$\beta = \alpha + 1$.  By the induction hypothesis, every abelian group in
$\PHKh_\alpha \frak F$ is of cardinality strictly less than
$\aleph_{\omega_0\alpha}$. If we apply Corollary~7.8, with
$\frak X =\PHKh_\alpha \frak F$, we see that every abelian group in
$\PHKh_1\frak X =\PHKh_\beta \frak F$ is of cardinality strictly less than
$\aleph_{\omega_0\alpha  + \omega_0} = \aleph_{\omega_0\beta}$.
\qed\enddemo

We can now make six statements, of varying profundity and novelty, about how
{\it free} abelian groups fit into the hierarchy $\PHKh\frak F$.

\proclaim
{7.10 Theorem} For each cardinal $\kappa$, let $\Bbb A_\kappa$ denote the free
abelian group of rank~$\kappa$.

\roster
\item $\Bbb A_0 \in \PHKh_0\frak{F}$.

\item $(\Bbb A_\kappa \mid 1 \le \kappa < \aleph_0)
\subseteq \PHKh_1\frak{F} - \PHKh_0\frak{F}$.

\item $(\Bbb A_\kappa \mid \aleph_0 \le \kappa < \aleph_{\omega_0})
\subseteq \PHKh_2\frak{F} - \PHKh_1\frak{F}$.

\item $\Bbb A_{\aleph_{\omega_0}}
\in \PHKh_3\frak{F} - \PHKh_2\frak{F}$.

\item For each finite ordinal $n$,
$\Bbb A_{\aleph_{\omega_0n}} \not \in \PHKh_{n+1}\frak F$; equivalently,
every free abelian group in $\PHKh_{n+1}\frak F$ has $\aleph$-rank strictly
less than
$\omega_0n$.

\item For each infinite ordinal $\beta$,
$\Bbb A_{\aleph_{\omega_0\beta}} \not \in \PHKh_{\beta}\frak F$.
\endroster
\endproclaim

\demo{Proof} (6) is a special case of Theorem~7.9.

(1) is clear.

(2).  Consider  $n \in \Bbb N$.  Then  $\Bbb A_n$ acts freely on an
$n$-dimensional CW-complex with underlying space $\Bbb R^n$, so
$\Bbb A_n \in \PHKh_1\frak F$.

It follows that $\Bbb A_n \in \PHKh_1\frak F$, and it is clear that, if $n
\ge 1$,
then  $\Bbb A_n \not\in \PHKh_0\frak F$.

(5). It is well known that, for each $n \in \Bbb N$,
$\H^n(\Bbb A_n, \Bbb Z)=\Bbb Z$.  This can be seen by induction, using a
long exact
sequence in cohomology which gives a short exact sequence of graded groups
$$0 \to \H^*(\Bbb A_{n-1}, \Bbb Z) \to \H^{*+1}(\Bbb A_{n}, \Bbb Z)
\to  \H^{*+1}(\Bbb A_{n-1}, \Bbb Z) \to 0;$$  it can also be seen from the
fact that
the $n$-torus is a $\operatorname{K}(\Bbb A_n,1)$.

Suppose that $\Bbb A_{\aleph_0} \in \PHKh_1\frak F$, so, for some $n \in
\Bbb N$,
$\Bbb A_{\aleph_0}$ acts freely on a contractible $n$-dimensional
CW-complex $X$.
Hence $\Bbb A_{n+1}$ acts freely on $X$, so $\H^{n+1}(\Bbb A_{n+1}, -) =
0$, which
contradicts the fact that $\H^{n+1}(\Bbb A_{n+1}, \Bbb Z) = \Bbb Z$.  Thus
$\Bbb A_{\aleph_0} \not\in \PHKh_1\frak F$, so every free abelian group in
$\PHKh_1\frak F$ has finite rank, that is, has $\arank$ equal to $-1$.

Now suppose that $n \ge 1$, and that every free abelian group in
$\PHKh_{n}\frak F$
has $\arank$ strictly less than $\omega_0(n-1)$.  Corollary~7.8, with
$\frak X = \PHKh_{n}\frak F$, and $\lambda = \omega_0(n-1)$, implies that
every free abelian group in $\PHKh_{n+1}\frak F$ has $\arank$ strictly less than
$\omega_0n$.

Now (5) has been proved by induction.

(3) and (4). By (2), $\PHKh_1\frak F$ includes all finitely generated free
abelian
groups.  Hence every free abelian group lies in $\PHKl\PHKh_1\frak F$.  If
we apply
Theorem~7.2(2) and (3), with $\frak X = \PHKh_1\frak F$, we find that the free
abelian groups of $\arank$ strictly less than $\omega_0$ lie in
$\PHKh_2\frak F$,
and that the free abelian group of $\arank$ $\omega_0$ lies in
$\PHKh_3\frak F$.  Combined with (5), these imply (3) and (4).
\qed\enddemo

\definition
{7.11 Remark} Theorem~7.10 gives an interesting new proof that
$$\PHKh_0\frak{F}\ne \PHKh_1\frak{F}\ne\PHKh_2\frak{F}\ne \PHKh_3\frak{F}.\qed$$
\enddefinition

\definition
{7.12 Conjecture}  $\Bbb A_{\aleph_{\omega_0+1}} \not \in \PHKh_3\frak{F}$;
equivalently, $\Bbb A_{\aleph_{\omega_0+1}} \not \in \PHKh\frak{F}$.
\qed \enddefinition

\definition
{7.13 Conjecture}  $\PHKh_3\frak{F} \ne \PHKh\frak{F}$.
\qed \enddefinition

\definition
{7.14 Conjecture}  There exists an ordinal $\alpha$ such that
$\PHKh_\alpha\frak{F} = \PHKh\frak{F}$.
\qed \enddefinition

\remark{\it Acknowledgments}  W\.~Dicks was partially supported by the DGES
and the DGI through grants PB96-1152 and BFM2000-0354, respectively.
I\.~Leary was
partially supported by the EPSRC.  S\.~Thomas was partially supported by
NSF grants.

I\.~Leary and P\.~Kropholler thank the Centre de Recerca Matem\`atica of the
Institut d'Estudis Catalans for the hospitality they received, and
S\.~Thomas thanks
the Departament de Matem\`atiques of the Universitat Aut\`onoma de
Barcelona for its
hospitality.
\endremark

\enddocument